\documentclass[11pt]{amsart}

\usepackage{amssymb,amsthm,amsmath}
\usepackage{color}
\definecolor{darkblue}{rgb}{0, 0, .4}
\definecolor{grey}{rgb}{.7, .7, .7}
\usepackage[breaklinks]{hyperref}
\hypersetup{
	colorlinks=true,
	linkcolor=darkblue,
	anchorcolor=darkblue,
	citecolor=darkblue,
	pagecolor=darkblue,
	urlcolor=darkblue,
	pdftitle={},
	pdfauthor={}
}
	
\newtheorem{theorem}{Theorem}[section]
\newtheorem{lemma}[theorem]{Lemma}

\theoremstyle{definition}
\newtheorem{definition}[theorem]{Definition}
\newtheorem{example}[theorem]{Example}

\theoremstyle{remark}

\numberwithin{equation}{section}

\theoremstyle{theorem}
\newtheorem{corollary}[theorem]{Corollary}

\newtheorem{proposition}[theorem]{Proposition}

\newtheorem{question}[theorem]{Question}

\newcommand{\N}[0]{\mathbb{N}}

\newcommand{\s}[0]{\sigma}
\newcommand{\g}[0]{\gamma}

\newcommand{\intersect}{\cap}

\newcommand{\gn}{ \bullet }  

\newcommand{\hs}{}  
\newcommand{\hd}{\diamond}  
\newcommand{\hz}{\circ}  
\newcommand{\hf}{\bullet}  
\newcommand{\hb}{{\color{grey} \bullet}}  
\newcommand{\hv}{\ast}  
\newcommand{\ha}{{\color{grey} \star}}  
\newcommand{\hh}{\star}  

\newcommand{\stringlessheap}{ \xymatrix @=-5pt @! } 
\newcommand{\heap}{ \xymatrix @=-9pt @! } 
\def\SDSize{6}  
\def\SDSizeTANRIGHT{4}  
\def\SDSizeTANLEFT{2}  
\def\SDMidpt{3}  
\def\SDColor{blue}
\def\SDEColor{black}
\newcommand{\StringRXD}[1]{{\color{\SDColor}\xy (\SDSize, \SDSize)*{}; (0, 0)*{}; **\crv{~**\dir{.}(\SDMidpt,\SDMidpt)};  (\SDMidpt, \SDMidpt)*{{\color{\SDEColor}#1}}; \endxy }}
\newcommand{\StringLXD}[1]{{\color{\SDColor}\xy (0, \SDSize)*{}; (\SDSize, 0)*{}; **\crv{~**\dir{.}(\SDMidpt,\SDMidpt)};  (\SDMidpt, \SDMidpt)*{{\color{\SDEColor}#1}}; \endxy }}
\newcommand{\StringLXDR}[1]{{\color{\SDColor}\xy (0, \SDSize)*{}; (\SDSize, 0)*{}; **\crv{~**\dir{.}(\SDMidpt,\SDMidpt)};  (\SDSize, \SDSize)*{}; (0, 0)*{}; **\crv{(\SDMidpt,\SDMidpt)}; (\SDMidpt, \SDMidpt)*{{\color{\SDEColor}#1}}; \endxy }}
\newcommand{\StringL}[1]{{\color{\SDColor}\xy (\SDSize, \SDSize)*{}; (0, 0)*{}; (0, \SDSize)*{}; **\crv{(\SDSizeTANLEFT,\SDMidpt)};  (\SDMidpt, \SDMidpt)*{{\color{\SDEColor}#1}}; \endxy}}
\newcommand{\StringR}[1]{{\color{\SDColor}\xy (0,0)*{}; (\SDSize, 0)*{}; (\SDSize, \SDSize)*{}; **\crv{(\SDSizeTANRIGHT,\SDMidpt)}; (\SDMidpt, \SDMidpt)*{{\color{\SDEColor}#1}}; \endxy }}
\newcommand{\StringLX}[1]{{\color{\SDColor}\xy (0, \SDSize)*{}; (\SDSize, 0)*{}; **\crv{(\SDMidpt,\SDMidpt)};  (\SDMidpt, \SDMidpt)*{{\color{\SDEColor}#1}}; \endxy }}
\newcommand{\StringRX}[1]{{\color{\SDColor}\xy (\SDSize, \SDSize)*{}; (0, 0)*{}; **\crv{(\SDMidpt,\SDMidpt)};  (\SDMidpt, \SDMidpt)*{{\color{\SDEColor}#1}}; \endxy }}
\newcommand{\StringLR}[1]{{\color{\SDColor}\xy (0, 0)*{}; (0, \SDSize)*{}; **\crv{(\SDSizeTANLEFT,\SDMidpt)};  (\SDSize, 0)*{}; (\SDSize, \SDSize)*{}; **\crv{( \SDSizeTANRIGHT,\SDMidpt)}; (\SDMidpt, \SDMidpt)*{{\color{\SDEColor}#1}}; \endxy }}
\newcommand{\StringLRX}[1]{{\color{\SDColor}\xy (0, \SDSize)*{}; (\SDSize, 0)*{}; **\crv{(\SDMidpt,\SDMidpt)};  (\SDSize, \SDSize)*{}; (0, 0)*{}; **\crv{(\SDMidpt,\SDMidpt)}; (\SDMidpt, \SDMidpt)*{{\color{\SDEColor}#1}}; \endxy }}

\usepackage[all, dvips, knot]{xy} 
\usepackage[all, knot, dvips, color]{xypic} 
\xyoption{arc}

\newcommand{\w}{\mathsf{w}}

\begin{document}

\title{Leading coefficients of Kazhdan--Lusztig polynomials for Deodhar elements}

\begin{abstract}
We show that the leading coefficient of the Kazhdan--Lusztig polynomial
$P_{x,w}(q)$ known as $\mu(x,w)$ is always either 0 or 1 when $w$ is a Deodhar
element of a finite Weyl group.  The Deodhar elements have previously been
characterized using pattern avoidance in \cite{b-w} and \cite{b-j1}.  In type
$A$, these elements are precisely the 321-hexagon avoiding permutations.  Using
Deodhar's algorithm \cite{d}, we provide some combinatorial criteria to determine
when $\mu(x,w) = 1$ for such permutations $w$.
\end{abstract}

\author{Brant C. Jones}
\address{Department of Mathematics, One Shields Avenue, University of California, Davis, CA 95616}
\email{\href{mailto:brant@math.ucdavis.edu}{\texttt{brant@math.ucdavis.edu}}}
\urladdr{\url{http://www.math.ucdavis.edu/\~brant/}}

\thanks{The author received support from NSF grants DMS-9983797 and DMS-0636297.}


\date{\today}

\maketitle


\bigskip
\section{Introduction to the 0-1 property}\label{s:mu.background}

In this paper we show that when $w$ is a Deodhar element of a finite Weyl group,
the coefficient of $q^{{1 \over 2}(l(w)-l(x)-1)}$ in the Kazhdan--Lusztig
polynomial $P_{x,w}(q)$ is always either 0 or 1.  This notorious coefficient is
known as $\mu(x,w)$ in the literature and corresponds to the term of highest
possible degree in $P_{x,w}(q)$.  The $\mu(x,w)$ values appear in multiplication
formulas for the Kazhdan--Lusztig basis elements $\{ {C'}_w \}$ of the Hecke
algebra and are used in the construction of Kazhdan--Lusztig graphs and cells
which in turn are used to construct Hecke algebra representations.  This was
first done in \cite{k-l} and is described in detail for type $A$ by \cite{b-b}.
Note that the $\mu$ coefficients control the recursive structure of the
Kazhdan--Lusztig polynomials and computing $\mu(x,w)$ is not known to be any
easier than computing the entire Kazhdan--Lusztig polynomial $P_{x,w}(q)$.  The
Deodhar elements have been characterized previously by pattern avoidance as
those that are [321]-hexagon avoiding permutations in \cite{b-w} and for the
other finite Weyl groups in \cite{b-j1}.

Our main result is motivated by the following open problem.

\begin{question}\label{q:main}
Fix a finite Coxeter group $W$.  Is there a simple characterization of the
elements $w \in W$ such that $\mu(x,w) \in \{0,1\}$ for all $x \in W$?
\end{question}

When all of the elements $w$ in a left cell $\mathcal{C}$ of the finite Weyl
group $W$ satisfy $\mu(x,w) \in \{0, 1\}$ for all $x \in W$, then
Kazhdan--Lusztig's construction of the representation of $W$ associated to
$\mathcal{C}$ depends only on the underlying graph structure of the $W$-graph,
rather than the edge-labeling of the graph by $\mu$ coefficients.  See
\cite{k-l} or Chapter 6 of \cite{b-b} for details.

Until fairly recently, it was conjectured that all $\mu(x,w)$ for type $A$ were
either 0 or 1.  See \cite{MW02} for some history about this conjecture and an
example in $S_{10}$ that shows the conjecture to be false.  On the other hand,
if $w$ and $x$ are elements of a finite Weyl group then under each of the
hypotheses below we have $\mu(x,w) \in \{0,1\}$.
\begin{enumerate}
\item[(1)]  Suppose $w \in S_n$, the symmetric group with $n \leq 9$.
\item[(2)]  Suppose $w$ is a Grassmannian permutation (equivalently, there is at most one decreasing consecutive pair of entries in the 1-line notation for $w$).
\item[(3)]  Suppose $w$ corresponds to a smooth Schubert variety.
\item[(4)]  Suppose $w$ and $x$ are fully commutative.
\item[(5)]  Suppose $w, x \in S_n$, the symmetric group, with $a(x) < a(w)$ where $a: S_n \rightarrow \N$ is the function defined by Lusztig in \cite{lusztig-a-function}.
\end{enumerate}
We have (1) from a difficult verification due to \cite{MW02} and (2) by
\cite{l-s}.  When the Schubert variety associated to $w$ is smooth, we have by
\cite{K-L2} that $P_{x,w}(q) = 1$ for all $x$, proving (3).  In type $A$, we
mention that $w$ is $\{[3412], [4231]\}$-avoiding if and only if $w$ is smooth
by \cite{la-s,ryan,wolper}.  In \cite{jgraham} (4) was observed for types $A$,
$D$ and $E$.  Recently, (4) has been extended to types $B$ and $H$ in
\cite{green4} and the proof techniques carry over to other Coxeter groups as
well.  We have (5) due to \cite{xi}.  Note that (4) and (5) require a
restriction on both $x$ and $w$.

Our main result in Theorem~\ref{t:mu.remain} adds a new set of elements to the
list above and will follow from Theorems~\ref{t:mu.main}, \ref{t:main.d} and
\ref{t:main.e} as well as Corollary~\ref{c:main.b} below.

\begin{theorem}\label{t:mu.remain}
If $w$ is a Deodhar element of a finite Weyl group $W$, and $x \in W$ then
$\mu(x,w) \in \{0, 1\}$.
\end{theorem}

In Section~\ref{s:background} we describe the main tools used in the proof of
Theorem~\ref{t:mu.remain}.  We will begin by proving Theorem~\ref{t:mu.remain}
for type $A$ in Sections~\ref{s:mu.a} and \ref{s:mu.a.proof}.  This proof also
holds for type $B$ since the unlabeled Coxeter graph of type $B$ is the same as
the unlabeled Coxeter graph of type $A$.  In Section~\ref{s:mu.d} we describe
how to reduce the type $D$ case to type $A$.  Finally, we verify the theorem
for the exceptional finite Weyl groups by computer in
Section~\ref{s:mu.others}.


\bigskip
\section{Background and introduction to $\mathbf{\mu}$-masks}\label{s:background}

We assume that the reader is familiar with the basic definitions and results
from Coxeter groups and Kazhdan--Lusztig theory discussed in \cite{h} or
\cite{b-b}, as well as the standard notions of pattern avoidance reviewed in
\cite{b-j1}.  Our main tools in this paper are Deodhar's combinatorial setting
\cite{d} for Kazhdan--Lusztig polynomials and the decorated heaps introduced in
\cite{b-w} and \cite{b-j1} to encode the masks of Deodhar's framework.  We refer
the reader to the introductions of \cite{b-w} and \cite{b-j1} for a more
leisurely discussion of these topics. 

\subsection{Deodhar's Theorem}

Let $w$ be an element of a Coxeter group $W$ and fix a reduced expression $\w =
\w_{1} \w_{2} \cdots \w_{k}$ for $w$, so each $\w_i$ is a Coxeter generator.
Define a \em mask \em $\s$ associated to the reduced expression $\w$ to be any
binary vector $(\s_1, \ldots, \s_k)$ of length $k = l(w)$.  Every mask
corresponds to a subexpression of $\w$ defined by $\w^\s = \w_{1}^{\s_1}
\cdots \w_{k}^{\s_k}$ where
\[
\w_{j}^{\s_j}  =
\begin{cases}
\w_{j}  &  \text{ if  }\s_j=1\\
\text{id}  &  \text{ if  }\s_j=0.
\end{cases}
\]
Each $\w^\s$ is a product of generators so it determines an element of $W$,
although $\w^\s$ may not be reduced.  For $1\leq j\leq k$, we also consider
initial sequences of masks denoted $\s[j] = (\s_1, \ldots, \s_j)$, and the
corresponding initial subexpressions $\w^{\s[j]} = \w_{1}^{\s_1} \cdots
\w_{j}^{\s_j}$.  In particular, we have $\w^{\s[k]} = \w^\s$.  The mask $\s$ is
\em proper \em if it does not consist of all 1 entries, since $\w^{(1, \ldots,
1)} = \w$ which is the fixed reduced expression for $w$.

We say that a position $j$ (for $2 \leq j \leq k$) of the fixed reduced
expression $\w$ is a \em defect \em with respect to the mask $\s$ if
\begin{equation*}
l(\w^{\s[j-1]} \w_{j}) < l(\w^{\s[j-1]}).
\end{equation*}
Note that the defect status of position $j$ does not depend on the value of
$\s_j$.  Let $d_{\w}(\s)$ denote the number of defects of $\w$ for a mask $\s$.
We will use the notation $d(\s) = d_{\w}(\s)$ if the reduced word $\w$ is fixed.  

Deodhar's framework gives a combinatorial interpretation for the
Kazhdan--Lusztig polynomial $P_{x,w}(q)$ as the generating function for masks
$\s$ on a reduced expression $\w$ with respect to the defect statistic $d(\s)$.
In this work, we consider the set
\[ \mathcal{S} = \{ 0, 1 \}^{l(w)} \]
of all $2^{l(w)}$ possible masks on $\w$ and define a prototype for
$P_{x,w}(q)$:
\[ P_x(q) = \sum_{ \substack{ \s \in \mathcal{S} \\ \w^{\s} = x } }
q^{d(\s)}. \]

\begin{definition}
Fix a reduced expression $\w$ for $w$.  We say that $\mathcal{S}$ is \em
bounded on \em $\w$ if $P_x(q)$ has degree $\leq {1 \over 2} (l(w) -
l(x) - 1)$ for all $x < w$ in Bruhat order.
\end{definition}

\begin{theorem}{\bf \cite{d}}\label{t:deodhar}
Let $x, w$ be elements in any Coxeter group $W$, and fix a reduced expression
$\w$ for $w$.  If $\mathcal{S}$ is bounded on $\w$ then 
\[ P_{x,w}(q) = P_x(q) = \sum_{ \substack{ \s \in \mathcal{S} \\ \w^{\s} = x } } q^{d(\s)}. \]
\end{theorem}

In \cite{b-w} and \cite{b-j1}, the elements $w$ of finite Weyl groups for which
$\mathcal{S}$ is bounded on any (equivalently, every) reduced expression $\w$ of
$w$ have been classified using pattern-avoidance.  We call such elements \em
Deodhar\em.  It follows from \cite[Corollary 5.3]{b-j1} that all of the Deodhar
elements are short-braid avoiding.  Here, we say that $w$ is \em short-braid
avoiding \em if there is no reduced expression for $w$ containing a factor of
the form $s t s$ where $s$ and $t$ are Coxeter generators satisfying $s t \neq t
s$.  We refer to the factor $s t s$ as a \em short-braid\em.  We define an
equivalence relation on the set of reduced expressions for a Coxeter element in
which two reduced expressions are in the same \em commutativity class \em if one
can be obtained from the other by a sequence of moves that interchange adjacent
Coxeter generators that commute.  Recall that by definition, every reduced
expression for a \em fully commutative \em element can be obtained from any
other using commuting moves on the Coxeter generators which appear in the
expressions.  In the simply-laced types $A$, $D$ and $E$, $w$ is fully
commutative if and only $w$ is short-braid avoiding.

By \cite[Lemma 2]{b-w}, we have that $\mathcal{S}$ is bounded on $\w$ if and
only if for every proper mask $\s \in \mathcal{S} \setminus \{ (1, 1, \ldots, 1)
\}$, we have
\begin{equation}\label{e:big_d}
D(\s) = \text{\# of plain-zeros of } \s - \text{\# of zero-defects of } \s > 0.
\end{equation}
We call $D(\s)$ the \em Deodhar statistic \em of $\s$.  Here, a position in $\w$
is a \em zero-defect \em if it has mask-value 0 and it is also a defect.  A
position in $\w$ is a \em plain-zero \em if it has mask-value 0 and it is not a
defect.  Similarly, we say that a position with mask-value 1 is a \em one-defect
\em if it is a defect, and we say that the position is a \em plain-one \em
otherwise.  

\subsection{Heaps}

If $\w = \w_1 \cdots \w_k$ is a reduced expression for a Coxeter element, then
following \cite{s1} we define a partial ordering on the indices $\{1, \ldots,
k\}$ by the transitive closure of the relation $i \lessdot j$ if $i < j$ and
$\w_i$ does not commute with $\w_j$.  We label each element $i$ of the poset by
the corresponding generator $\w_i$.  It follows quickly from the definition
that if $\w$ and $\w'$ are two reduced expressions for an element $w$ that are
in the same commutativity class then the labeled posets of $\w$ and $\w'$ are
isomorphic.  This isomorphism class of labeled posets is called the \em heap
\em of $\w$, where $\w$ is a reduced expression representative for a
commutativity class of $w$.  In particular, if $w$ is fully commutative then it
has a single commutativity class, and so there is a unique heap of $w$.

Cartier and Foata \cite{cartier-foata} were among the first to study heaps of
dimers, which were generalized to other settings by Viennot \cite{viennot}.
Stembridge has studied enumerative aspects of heaps \cite{s1,s2} in the context
of fully commutative elements.  Green has also considered heaps of pieces with
applications to Coxeter groups in \cite{green1,green2}.

The Coxeter graph of type $A_{n}$ has the form
\[
\xymatrix @-1pc {
\gn_1 & \ar@{-}[l] \gn_2 \ar@{-}[r] & \gn_3 \ar@{-}[r] & \ldots - \gn_{n}\\
} .
\]
As in \cite{b-w}, we will represent a heap of a type $A$ element as
a set of lattice points embedded in $\N^2$.  To do so, we assign coordinates
$(x,y) \in \N^2$ to each entry of the labeled Hasse diagram for the heap of $\w$
in such a way that:
\begin{enumerate} 
\item[(1)] An entry represented by $(x,y)$ is labeled by the Coxeter
generator $s_i$ in the heap if and only if $x = i$, and
\item[(2)] An entry represented by $(x,y)$ is greater than an entry
represented by $(x',y')$ in the heap if and only if $y > y'$.
\end{enumerate}
Since the Coxeter graph of type $A$ is a path, it follows from the definition
that $(x,y)$ covers $(x',y')$ in the heap if and only if $x = x' \pm 1$, $y >
y'$, and there are no entries $(x'', y'')$ such that $x'' \in \{x, x'\}$ and $y'
< y'' < y$.  Hence, we can completely reconstruct the edges of the Hasse
diagram and the corresponding heap poset from a lattice point representation.
This representation enables us to make arguments ``by picture'' that would
otherwise be difficult to formulate.

Although there are many coordinate assignments for any particular heap, the $x$
coordinates of each entry are fixed for all of them, and the coordinate
assignments of any two entries only differs in the amount of vertical space
between them.  In the case that $w$ is fully commutative, a canonical choice can
be made by ``coalescing'' the entries as in \cite{b-w}.

To carry this out, we form the heap of a reduced expression $\w$ by reading $\w$
from left to right, and dropping a point into the column representing each
generator $\w_i$.  We envision each point as being slightly wider than its
column and under the influence of gravity, in the sense that the point must fall
to the lowest possible position in the column over the generator corresponding
to $\w_i$ without passing any previously placed points in adjacent columns.  To
coalesce these points, we apply gravity in the other direction.  We say that two
points $(x,y)$, $(x',y')$ are \em connected \em if $x = x' \pm 1$ and $y = y'
\pm 1$.  Whenever there exists a connected component lying below another
connected component with empty lattice points between them, then we allow the
first component to rise up until it is blocked by the second component.  We
apply these elevations until the heap is pushed together as much as possible.  

In the example below, we show the heap of the reduced expression $\w = s_1 s_4
s_2 s_3 s_5$ before and after coalescing.
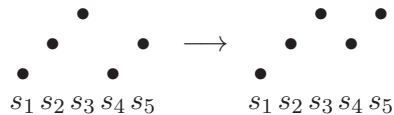
\[
\stringlessheap {
\hs &   & \hf & & \hs \\
    & \hf &   & \hs & \hf \\
\hf &   & \hs & \hf & \hs \\
s_1   & s_2  & s_3  & s_4 & s_5 & \hs \\ 
}
\parbox[t]{0.3in}{ \vspace{0.1in} $\longrightarrow$ }
\stringlessheap {
\hs &   & \hf & & \hf \\
    & \hf &   & \hf \\
\hf &   & \hs & \hs \\
s_1   & s_2  & s_3  & s_4 & s_5 \\ 
}
\]

In type $A$, the heap construction can be combined with another combinatorial
model for permutations in which the entries from the 1-line notation are
represented by strings.  Here, the \em 1-line notation \em $w = [w_1 w_2 \cdots
w_n]$ specifies the permutation $w$ as the bijection mapping $i$ to $w_i$.  The
points at which two strings cross can be viewed as adjacent transpositions of
the 1-line notation.  Hence, we may overlay strings on top of a heap diagram to
recover the 1-line notation for the element, by drawing the strings from bottom
to top so that they cross at each entry in the heap where they meet and bounce
at each lattice point not in the heap.  Conversely, each permutation string
diagram corresponds with a heap by taking all of the points where the strings
cross as the entries of the heap.

For example, we can overlay strings on the two heaps of the permutation $[3214]$
and the heap of $s_1 s_4 s_2 s_3 s_5$.  Note that the labels in the picture
below refer to the strings, not the generators.
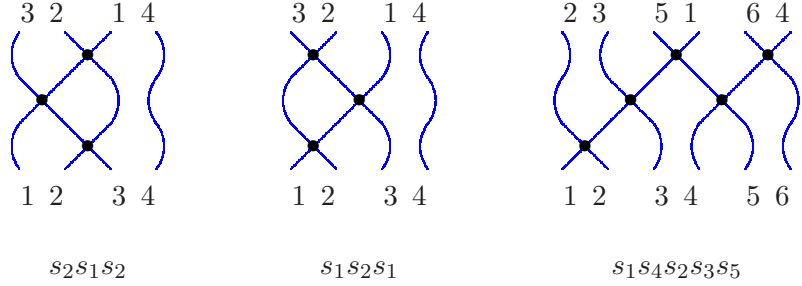
\begin{center}
\begin{tabular}{ccc}
	\heap{
	& 3 \ \ 2 &  & 1 \ \ 4 & \\
	\StringR{\hs} & \hs & \StringLRX{\hf} & \hs & \StringL{\hs} \\
	\hs & \StringLRX{\hf} & \hs & \StringLR{\hs} & \hs \\
	\StringR{\hs} & \hs & \StringLRX{\hf} & \hs & \StringL{\hs} \\
	& 1 \ \ 2 &  & 3 \ \ 4 & \\
	& & & & \\
	} & 
	\heap{
	& 3 \ \ 2 &  & 1 \ \ 4 &  \\
	\hs & \StringLRX{\hf} & \hs & \StringLR{\hs} \\
	\StringR{\hs} & \hs & \StringLRX{\hf} & \hs & \StringL{\hs} \\
	\hs & \StringLRX{\hf} & \hs & \StringLR{\hs} \\
	& 1 \ \ 2 &  & 3 \ \ 4 &  \\
	& & & & \\
	} & 
	\heap{
	& 2 \ \ 3 &  & 5 \ \ 1 & & 6 \ \ 4 \\
	& \StringLR{\hs} & & \StringLRX{\hf} & & \StringLRX{\hf} \\
	\StringR{\hs} & & \StringLRX{\hf} & & \StringLRX{\hf} & & \StringL{\hs} \\
	& \StringLRX{\hf} & & \StringLR{\hs} & & \StringLR{\hs} \\
	& 1 \ \ 2 &  & 3 \ \ 4 &  & 5 \ \ 6  \\
	& & & & \\
	} \\
$s_2 s_1 s_2$ & $s_1 s_2 s_1$ & $s_1 s_4 s_2 s_3 s_5$ \\
\end{tabular}
\end{center}

\subsection{Heaps of types $B$ and $D$}\label{s:type_bd}

The Coxeter graph of type $B_{n}$ is of the form
\[
\xymatrix @-1pc {
\gn_0 & \ar@{-}[l]_{4} \gn_1 \ar@{-}[r] & \gn_2 \ar@{-}[r] & \gn_3 \ar@{-}[r] 
& \ldots - \gn_{n-1}\\
} .
\]
From this graph, we see that the symmetric group $S_{n-1}$ is a parabolic subgroup
of this Coxeter group.  Because of this, the elements of this group have a
standard 1-line notation in which a subset of the entries are \textit{barred}.
We often think of the barred entries as negative numbers, and this group is
referred to as the group of \textit{signed permutations} or the
\textit{hyperoctahedral group}.  The action of the generators on the 1-line
notation is the same for $\{s_1, s_2, \ldots, s_{n-1} \}$ as in type $A$ in which
$s_i$ interchanges the entries in positions $i$ and $i+1$ in the 1-line
notation for $w$.  The $s_0$ generator acts on the right of $w$ by changing the
sign of the first entry in the 1-line notation for $w$.  For example,
$w=[\bar{4}2\bar{3}1]$ is an element of $B_{4}$ and
\begin{align*}
s_{0}&=[\bar{1}234]\\
ws_{0}&=[42\bar{3}1]\\
ws_{1}&=[2\bar{4}\bar{3}1].
\end{align*}

The Coxeter graph for type $D_{n}$ is
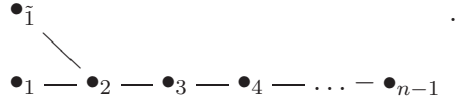
\[
\xymatrix @-1pc {
{\gn}_{\tilde{1}} &        &                &                            &  \\
\gn_1 & \ar@{-}[l] \ar@{-}[ul] \gn_2 \ar@{-}[r] & \gn_3 \ar@{-}[r] & \gn_4 \ar@{-}[r] & \ldots - \gn_{n-1} \\
} .
\]
The elements of type $D$ can be viewed as the subgroup of $B_{n}$ consisting of
signed permutations with an even number of barred entries.  The action of the
generators on the 1-line notation is the same for $\{s_1, s_2, \ldots \}$ as in
type $A$ in which $s_i$ interchanges the entries in positions $i$ and $i+1$ in
the 1-line notation for $w$.  The $s_{\tilde{1}}$ generator acts on the right of
$w$ by marking the first two entries in the 1-line notation for $w$ with bars
and interchanging them.  For example, $w=[\bar{4}2\bar{3}1]$ is an element of
$D_{4}$ and
\begin{align*}
w s_{\tilde{1}} =& [\bar{2}4\bar{3}1] .
\end{align*}
Although the Coxeter graph for type $D$ has a fork, we will assign coordinates
to draw the heap of type $D$ elements in a linearized way by allowing entries in
the first column to consist of either generator $s_1$, $s_{\tilde{1}}$, or both:
\begin{align*}
s_{1}=& \hf\\
s_{\tilde{1}}=& \tilde{\hf} \\
s_{1}s_{\tilde{1}}=& \hf \tilde{\hf} .
\end{align*}
Hence, if $w \in D_n$ then the heap of $\w$ has coordinates from $[n-1] \times
\N$, where $[n-1] = \{1, 2, \ldots, n-1\}$ and lattice points in the first
column may be assigned to two entries of the heap.

As in type $A$, we can adorn the heap of $\w$ with strings that represent the
entries of the 1-line notation for the element.  If we label the strings at the
bottom of the diagram with the numbers from $1, \ldots, n$ then the
$s_{\tilde{1}}$ generator crosses the strings which intersect it and changes
the sign on the labels for both strings.  All other generators simply cross the
strings as in type $A$.  For example, the heap of the reduced expression
$s_{\tilde{1}}s_{2}s_{3}s_{1}s_{2}s_{\tilde{1}}s_{1}$
is shown below.
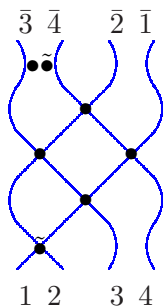
\[
	\heap {
	& \bar{3} \ \ \bar{4} & & \bar{2} \ \ \bar{1} &  &  \\
	& \StringLR{\hf \tilde{\hf}} & \hs & \StringLR{\hs} & \hs & \hs & \hs \\
	\StringR{\hs} & & \StringLRX{\hf} & \hs & \StringL{\hs} & \hs & \hs & \hs \\
	& \StringLRX{\hf} & {\hs} & \StringLRX{\hf} & \hs & \hs & \hs \\
	\StringR{\hs} & & \StringLRX{\hf} & \hs & \StringL{\hs} & \hs & \hs & \hs \\
	& \StringLRX{\tilde{\hf}} & \hs & \StringLR{\hs} & \hs & \hs & \hs \\
	& 1 \ \ 2 & & 3 \ \ 4  \\
	}
\]

Let $\w$ be a reduced expression for an element $w$ of type $A$ or type $D$.
Suppose $x$ and $y$ are a pair of entries in the heap of $\w$ that correspond to
the same generator $s_i$, so they lie in the same column $i$ of the heap.
Assume that $x$ and $y$ are a \em minimal pair \em in the sense that there is no
other entry between them in column $i$.  Then, for $\w$ to be reduced, there
must exist at least one generator that does not commute with $s_i$ lying between
$x$ and $y$, and if $\w$ is short-braid avoiding, there must actually be two
entries that do not commute with $s_i$ and lie strictly between $x$ and $y$ in
the heap.  We call these two non-commuting labeled heap entries a \em resolution
\em of the pair $x,y$.  If both of the generators in a resolution lie in column
$i-1$ ($i+1$, respectively), we call the resolution a \em left (right,
respectively) resolution\em.  If the generators lie in distinct columns, we call
the resolution a \em distinct resolution\em.  The Lateral Convexity Lemma
\cite[Lemma 1]{b-w} characterizes fully commutative permutations $w$ as those
for which every minimal pair in the heap of $w$ has a distinct resolution.  If
$w \in D_{n}$ is short braid avoiding and every minimal pair of entries in the
heap of $w$ has a distinct resolution then we say $w$ is \em convex\em.  The
heap of a convex element in type $D$ can be coalesced just as in type $A$; see
\cite[Section 7]{b-j1} for details.

We will frequently interpret a mask applied to a reduced expression for a
Deodhar element as a decorated heap diagram, and overlay strings on this
decorated heap.  In particular, since Deodhar elements are fully commutative, we
may refer to a mask $\s$ on the element $w$ and draw the unique decorated heap
associated to $\s$, rather than explicitly fixing a reduced expression $\w$ for
each mask we consider.  If we allow the strings to cross exactly at mask-value 1
entries of the decorated heap then the strings encode the 1-line notation for
$\w^{\s}$.  In the figures that follow, we decorate the heap diagrams according
to mask-value using Table~\ref{ta:1}.

\bigskip
\begin{center}
\begin{table}[ht]
\centering
\caption{Heap decorations}
\begin{tabular}{|p{0.7in}|p{3in}|} 
\hline
Decoration & Mask-value \\
\hline
$\hd$  &	zero-defect entry\\
$\hz$  &	plain-zero entry (not a defect)\\
$\hf$  &	mask-value 1 entry\\
$\ha$  &	entry of the heap with unknown mask-value\\
$\hb$  &	lattice point not necessarily in the heap\\
$\hv$  &	lattice point that is definitely not in the heap\\
$\hh$  &	lattice point that is highlighted for emphasis\\
\hline
\end{tabular}
\label{ta:1}
\end{table}
\end{center}
\bigskip

In type $A$, a decorated heap entry is a defect precisely when the two strings
emanating from the entry cross below the entry an odd number of times.  This
follows because an odd string crossing below the entry $\w_j$ in the decorated
heap of $\w$ corresponds to a descent in the 1-line notation for $\w^{\s[j-1]}$.
Hence, each defect $d$ must have a \em left \em and \em right critical zero \em
along the paths of the left and right strings respectively, which allow the
strings to eventually cross.  These critical zeros are denoted $lcz(d)$ and
$rcz(d)$ respectively, and when $d$ has mask-value 0 then we declare that $d$ is
a critical zero of itself as well.  For example, we have $lcz(d) = p$ and
$rcz(d) = q$ in the heap shown below.
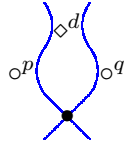
\[ \heap {
& {\hs} & & \StringLR{\hd^d} & & {\hs} \\
{\hs} & & \StringR{\hz^p} & & \StringL{\hz^q} & & {\hs} \\
& {\hs} & & \StringLRX{\hf} & & {\hs} & & {\hs} \\
} \]
The analogous result using signed strings to determine the defect status of a
type $D$ heap entry is given in \cite[Lemma 6.1]{b-j1}.  When the entry being
tested does not correspond to $s_{\tilde{1}}$, the test is simply the signed
version of the type $A$ inversion test.

Recall the following heap fragments called the $I$-\textit{shape} and the
4-\textit{stack} indicated with black entries below.  These fragments were
introduced for type $D$ in \cite{b-j1}.
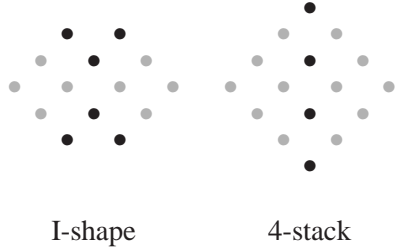
\begin{equation*}
\begin{matrix}
        \xymatrix @=-2pt @! {
& {\hs} & & {\hs} & & {\hs} \\
{\hs} & & {\hf} & & {\hf} & & {\hs} \\
& {\hb} & & {\hf} & & {\hb}\\
{\hb} & & {\hb} & & {\hb} & & {\hb} \\
& {\hb} & & {\hf} & & {\hb} \\
{\hs} & & {\hf} & & {\hf} & & {\hs} \\
& {\hs} & & {\hs} & & {\hs} \\
    }  &
         \xymatrix @=-2pt @! {
& {\hs} & & {\hf} & & {\hs} \\
{\hs} & & {\hb} & & {\hb} & & {\hs} \\
& {\hb} & & {\hf} & & {\hb}\\
{\hb} & & {\hb} & & {\hb} & & {\hb} \\
& {\hb} & & {\hf} & & {\hb} \\
{\hs} & & {\hb} & & {\hb} & & {\hs} \\
& {\hs} & & {\hf} & & {\hs} \\
    }  \\
    & \\
    \text{I-shape} & \text{4-stack} \\
\end{matrix}
\end{equation*}

\begin{lemma}\label{l:shape}
If $w$ is a Deodhar element in type $A$ or a convex Deodhar element of type $D$
with entries from at least two levels in the first column of the heap of $w$,
then the heap of $w$ does not contain an I-shape or a 4-stack.
\end{lemma}
\begin{proof}
By lateral convexity \cite[Lemma 1]{b-w}, any type $A$ heap which contains
entries in either of the configurations must contain a hexagon and is
consequently not Deodhar.  This result is extended in \cite[Lemma 8.5]{b-j1} to
convex elements of type $D$ with entries from at least two levels in the first
column.
\end{proof}

We say that a 3-\textit{stack} is the heap fragment obtained from the 4-stack by
removing the top or bottom entry.

\subsection{$\mathbf{\mu}$-masks for elements of finite Weyl groups}

We begin by interpreting $\mu(x,w)$ in Deodhar's framework.

\begin{definition}\label{d:mu_mask}
Fix a reduced expression $\w$ for a Deodhar element $w$ in any Coxeter group.
We say that a proper mask $\s \in \mathcal{S}$ is a \em $\mu$-mask for $\w$ \em
if
\begin{equation}
\text{\# defects of } \s = d(\s) = \frac{l(\w)-l(\w^{\s})-1}{2} .
\end{equation}
\end{definition}

A term contributing to the $\mu$ coefficient occurs in degree ${1 \over 2}
(l(\w) - l(\w^{\s}) - 1)$ of the Kazhdan--Lusztig polynomial $P_{\w^{\s},
\w}(q)$.  Hence, by Theorem~\ref{t:deodhar} we have 
\[ \mu(x, w) = \#\{ \s \in \mathcal{S} | \w^{\s} = x \text{ and } \s \text{ is a
$\mu$-mask for } \w \}.  \]

\begin{lemma}\label{l:mu_mask}
Let $\w$ be a reduced expression for a Deodhar element $w$ in any Coxeter group.
Then, $\s$ is a $\mu$-mask for $\w$ if and only if 
\[ \text{\# of plain-zeros of } \s - \text{\# of zero-defects of } \s = D(\s) = 1. \]
\end{lemma}

The following proof uses a similar technique to \cite[Lemma 2]{b-w}. 

\begin{proof}
Observe that
\begin{align*}
& {1 \over 2} (l(\w) - l(\w^{\s}) - 1) = d(\s) \\
\iff \ \ \ \ & l(\w) - l(\w^{\s}) - 2 d(\s) = 1.
\end{align*}
We show that
\begin{equation}\label{e:defect_translation}
l(\w) - l(\w^{\s}) - 2 d(\s) = D(\s)
\end{equation}
by induction on the length of $\w$.  

One can check that \eqref{e:defect_translation} holds when $l(\w) = 1$.  Suppose
\eqref{e:defect_translation} holds for all $\w$ such that $l(\w) \leq k$.  Next,
let $\w$ be a reduced expression for a Deodhar element of length $k$ and
suppose that $\w s_i > \w$, so $\w s_i$ is a reduced expression of length $k+1$.
Every mask for $\w s_i$ is obtained from a mask on $\w$ extended by 0 or 1 at
the last generator.  Consider the effect on the left and right sides of
\eqref{e:defect_translation} respectively, in each of the possible cases as
shown in Table~\ref{ta:2}.

\medskip
\begin{table}[h]
\centering
\caption{Multiplication by $s_i$}
\begin{tabular}{|llll|}
\hline
Mask-value of $s_i$ & Defect-status of $s_i$ & Change in left side & Change in
right side \\
\hline
1 & not a defect & $1-1 = 0$ & $0$ \\
1 & defect & $1 - (-1) - 2 (1) = 0$ & $0$ \\
0 & not a defect & $1$ & $1$ \\
0 & defect & $1 - 2 (-1) = -1$ & $-1$ \\
\hline
\end{tabular}
\label{ta:2}
\end{table}
\medskip

By induction, the left and right sides of \eqref{e:defect_translation} are
equal for arbitrary $\w$ and $\s$.
\end{proof}

\begin{lemma}\label{l:no_isolation}
Let $\w$ be a Deodhar element in any Coxeter group and $\s$ be a $\mu$-mask on
$\w$.  If $\s'$ is any mask obtained from $\s$ by changing the mask-value of a
plain-zero to have mask-value 1, then $\s'$ must have strictly fewer
zero-defects than $\s$.  Also, there are no maximal entries in the heap of $\w$
that are one-defects in $\s$.
\end{lemma}
\begin{proof}
If $\s'$ does not have fewer zero-defects than $\s$ then we have
\[
D(\s') = \text{ \# plain-zeros } \s' - \text{ \# zero-defects } \s'
< \text{ \# plain-zeros } \s - \text{ \# zero-defects } \s = D(\s) = 1
\]
so $\s'$ is a non-Deodhar mask for $\w$ by Equation~(\ref{e:big_d}), which is a
contradiction.

Similarly, if there is a maximal one-defect $d$ in the decorated heap of $\s$
then we can form a non-Deodhar mask $\s''$ on $\w$ by changing the mask-value
of $d$ to 0.  To see this, observe that $d$ is a zero-defect in $\s''$ and
every other entry in $\s''$ has the same defect status as in $\s$ since $d$ is
heap-maximal.  This implies that $D(\s'') < D(\s) = 1$ because $\s$ is a $\mu$-mask.
\end{proof}

\begin{lemma}\label{l:mu.winverse}
Let $W$ be a Coxeter group.  For any $w, x \in W$ we have $\mu(x,w) =
\mu(x^{-1}, w^{-1})$.
\end{lemma}
\begin{proof}
This follows from the standard fact that $P_{x,w}(q) = P_{x^{-1}, w^{-1}}(q)$;
see \cite[Chapter 5]{b-b} for example.
\end{proof}


\bigskip
\section{Properties of $\mu$-masks in type $A$}\label{s:mu.a}

In this section we restrict to type $A$.  Recall from \cite{b-w} the construction
of the \em defect graph \em $G_{\s}$ associated to $\s$.  The vertices of
$G_{\s}$ consist of the zero-defects of $\s$, and there exists an edge between
two entries $x$ and $y$ if they share any critical zeros:
\[ \{ lcz(x), rcz(x), x \} \intersect \{ rcz(y), lcz(y), y \} \neq \emptyset. \]

\begin{example}\label{e:typical}
Here is a decorated heap corresponding to a typical $\mu$-mask for a Deodhar
permutation.
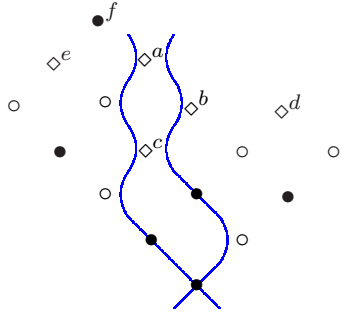
\[ \heap {
{\hs} & & {\hf^f} & & {\hs} & & {\hs} \\
& {\hd^e} & & \StringLR{\hd^a} & & {\hs} \\
{\hz} & & \StringR{\hz} & & \StringL{\hd^b} & & {\hd^d} \\
& {\hf} & & \StringLR{\hd^c} & & {\hz} & & {\hz} \\
{\hs} & & \StringR{\hz} & & \StringLX{\hf} & & {\hf} \\
& {\hs} & & \StringLX{\hf} & & \StringL{\hz} \\
{\hs} & & {\hs} & & \StringLRX{\hf} & & {\hs} \\
} \]
The defect graph for this mask is
\[ \xymatrix @-1pc {
\gn_e \ar@{-}[r] & \gn_a \ar@{-}[r] & \gn_b \ar@{-}[d] \ar@{-}[r] & \gn_d \\ 
 &  & \gn_c &  \\ 
} . \]
\end{example}

\begin{lemma}{\bf \cite[Proposition 1, Lemma 5]{b-w}}\label{l:b-w}
Suppose $w$ is a Deodhar permutation and $\s$ is any mask on $w$.  Then,
\begin{enumerate}
\item[(1)]  The defect graph $G_{\s}$ is a forest.
\item[(2)]  In the decorated heap of $\s$, no entry is a critical zero for more
	than two zero-defects, including itself.
\end{enumerate}
\end{lemma}

\begin{figure}[ht]
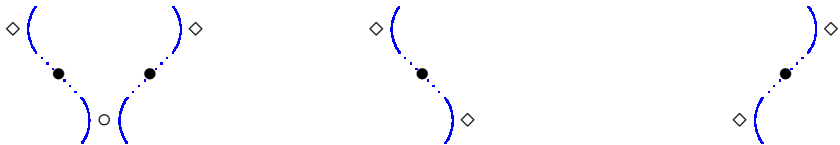

\begin{center}
$$\begin{matrix} 
\heap {
{\hs} & & \StringR{\hd} & & {\hs} & & \StringL{\hd} \\
& {\hs} & & \StringLXD{\hf} & & \StringRXD{\hf} & & {\hs} \\
{\hs} & & {\hs} & & \StringLR{\hz} & & {\hs} \\
} &
\heap {
{\hs} & & \StringR{\hd} & & {\hs} & & {\hs} \\
& {\hs} & & \StringLXD{\hf} & & {\hs} & & {\hs} \\
{\hs} & & {\hs} & & \StringL{\hd} & & {\hs} \\
} & 
\heap {
{\hs} & & {\hs} & & \StringL{\hd} & & {\hs} \\
& {\hs} & & \StringRXD{\hf} & & {\hs} & & {\hs} \\
{\hs} & & \StringR{\hd} & & {\hs} & & {\hs} \\
} \\
\end{matrix}$$
\end{center}
\caption{Zero-defect connections}
\label{f:zd_connections}
\end{figure}

Lemma~\ref{l:b-w}(2) implies that when $\s$ is a mask on a Deodhar permutation,
the only possible connections between zero-defects in $G_{\s}$ are those three
shown in Figure~\ref{f:zd_connections}.

The following result describes the defect graphs of $\mu$-masks.

\begin{lemma}\label{l:pc}
Let $\s$ be a mask for a Deodhar permutation $w$.  We have that $\s$ is a
$\mu$-mask if and only if $G_{\s}$ is a tree and every plain-zero in $\s$ is a
critical zero for some zero-defect.
\end{lemma}
\begin{proof}
Suppose $\s$ is a mask for $w$.  Every zero-defect has 3 critical zeros
including itself, and no point is a critical zero for more than two
zero-defects by Lemma~\ref{l:b-w}(2).  Hence, the number of critical zeros
that are shared by some pair of zero-defects is given by the number of edges in
$G_{\s}$.  Let $V$ and $E$ be the number of vertices and edges in $G_{\s}$
respectively.  Let $X$ be the number of plain-zeros that are not critical
for any zero-defect in $\s$.  Then we have
\begin{align*}
\text{ \# plain-zeros of } \s + \text{ \# zero-defects of } \s & = 3 V - E + X \\
\text{ \# zero-defects in } \s & = V \\
\end{align*}
so,
\begin{align*}
\text{ \# plain-zeros in } \s & = 2 V - E + X \\
D(\s) = \text{ \# plain-zeros} - \text{\# zero-defects in } \s & = V - E + X . \\
\end{align*}
If $\s$ is a $\mu$-mask then by Lemma~\ref{l:mu_mask} we have
\[ V - E + X = 1 . \]
By Lemma~\ref{l:b-w}(1) we have that $G_{\s}$ is a forest, so $E \leq V
- 1$.  This implies $X \leq 0$ so $X = 0$.  Furthermore, $E = V - 1$ so
$G_{\s}$ is connected.

On the other hand, if $G_{\s}$ is a tree then $V - E = 1$.  Also, if every plain-zero
is a critical zero for some zero-defect then $X = 0$.  Hence, we have
\[ D(\s) = V - E  + X = 1 \]
so $\s$ is a $\mu$-mask by Lemma~\ref{l:mu_mask}.
\end{proof}

In the proofs that follow, we will frequently use the fact from
Lemma~\ref{l:pc} that $G_{\s}$ is connected if $\s$ is a $\mu$-mask to
contradict the hypothesis that a given mask is actually a $\mu$-mask.  We use
the following definition to give a local reformulation of this result.  

\begin{definition}\label{d:sep_trip}
In a decorated coalesced heap fragment for a Deodhar permutation, we say that a
connected collection of entries lying on a diagonal in the lattice $\N^2$ is a
\em diagonal \em of the heap fragment.

Let $\s$ be a mask on a Deodhar permutation $w$, and suppose there exists a set
of entries $H$ from three parallel adjacent diagonals in the heap of $w$ whose
mask-value/defect-status are known.  Moreover, suppose that these diagonals all
pass through a common column, so the heap of $w$ contains a 3-stack.  Also,
suppose that no zero-defect lying to the left of $H$ can share a critical zero
with any zero-defect lying to the right of $H$ without introducing another
diagonal that passes through the common column, in which case the heap contains
a 4-stack.  We call such a subset $H$ a \em separating triple of diagonals\em.  
\end{definition}

In practice, we find separating triples of diagonals by considering paths of
strings through the heaps as in Example~\ref{ex:sep}.

\begin{example}\label{ex:sep}
Consider the decorated heap fragment associated to $\s$ shown below.  Suppose
that $g$ is not in the heap, so $e$ cannot be connected to any zero-defects to
the right of $e$ in $G_{\s}$ by lateral convexity.  Also, suppose that $h$ is
not in the heap, so the left-critical zero of $d$ cannot be shared with a
zero-defect on the left.  Then, the three diagonal segments running in the
northwest-southeast direction form a separating triple of diagonals.

As we will see below, this implies that $\s$ is not a
$\mu$-mask, because there is no way to connect $d$ and $e$ in $G_{\s}$
regardless of the mask-value/defect-status of the other entries in $\s$ that are
not shown.
$$\begin{matrix} 
\heap {
{\hs} & & \StringLR{\hd^e} & & {\hv^g} & & {\hs} \\
& \StringR{\hz} & & \StringL{\hz} & & {\hs} & & {\hs} \\
{\hs} & & \StringLRX{\hf} & & {\hs} & & {\hs} \\
& {\hs} & & \StringLR{\hd^d} & & {\hs} & & {\hs} \\
{\hs} & & \StringR{\hz} & & \StringL{\hz} & & {\hs} \\
& {\hv^h} & & \StringLRX{\hf} & & {\hs} & & {\hs} \\
}\end{matrix}$$
\end{example}

Motivated by Example~\ref{ex:sep}, we make the following definition.

\begin{definition}\label{d:sep_pair}
Suppose there exists a separating triple of diagonals $H$ containing two
zero-defects $d$ and $e$ that are not connected in the induced defect graph of
$\s$ restricted to the entries from $H$.  If:
\begin{enumerate}
\item[(1)]  $d$ does not share any critical zero (including itself)
with a zero-defect lying to the left (or right) of $d$, and 
\item[(2)]  $e$ does not share any critical zero with a zero-defect lying to the
right (left, respectively) of $e$, 
\end{enumerate}
then we say that $d$ and $e$ form a \em separating pair in $G_{\s}$\em. 
\end{definition}

\begin{lemma}\label{l:connected}
Suppose $\s$ is a mask for a Deodhar permutation $w$.  If $d$ and $e$ are
entries of the decorated heap of $\s$ that form a separating pair in $G_{\s}$
then $\s$ is not a $\mu$-mask.
\end{lemma}
\begin{proof}
Assume $\s$ is a $\mu$-mask, so $d$ and $e$ are connected in $G_{\s}$ by
Lemma~\ref{l:pc}.  Suppose $d$ and $e$ are not connected in the defect graph of
$\s$ restricted to the entries $H$ that lie along the three diagonals whose
defect-status and mask-value are known.  Then there exist zero-defects $d'$ and
$e'$ outside of $H$ such that $d'$ is connected to $d$ in $G_{\s}$, $e'$ is
connected to $e$ in $G_{\s}$, and $d'$ and $e'$ are adjacent in $G_{\s}$ via
one of the connections shown in Figure~\ref{f:zd_connections}.  By the
assumptions (1) and (2) in Definition~\ref{d:sep_pair} on $d$ and $e$, we have
that $d'$ lies on the opposite side of $H$ from $e'$.  Hence, we have by
Definition~\ref{d:sep_trip} that $d'$ and $e'$ cannot share a critical zero
without introducing a 4-stack.  Since a 4-stack implies that $w$ is not Deodhar,
we have a contradiction.  Thus, the existence of a separating pair proves that
$\s$ is not a $\mu$-mask.
\end{proof}


\bigskip
\section{Proof of the 0-1 property for type $A$ Deodhar elements}\label{s:mu.a.proof}

We are now in a position to prove our main result for type $A$.

\begin{theorem}\label{t:mu.main}
Let $w$ be a Deodhar permutation and suppose $x$ is another permutation.  Then,
$\mu(x,w) \in \{0, 1\}$.
\end{theorem}
\begin{proof}
Let $w$ be a minimal length element such that $w$ is Deodhar and $\mu(x,w) > 1$
for some $x$.  Then there exist distinct $\mu$-masks $\s$ and
$\g$ such that $w^{\s} = x = w^{\g}$.  

Consider the coalesced heap of $w$.  By lateral convexity \cite[Lemma 1]{b-w},
$w$ has a single entry in the leftmost occupied column of the heap and since $w$
is Deodhar there are fewer than 4 entries in the fourth occupied column by
Lemma~\ref{l:shape}.  Hence, at least one of the points labeled $y$ or $y'$ is
not in the heap:
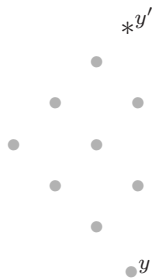
\[
\xymatrix @=-4pt @! {
& {\hs} & & {\hv}^{y'} & & {\hs} & & {\hs} \\
& & {\hb} & & {\hs} & & {\hs} & \\
& {\hb} & & {\hb} & & {\hs} & & {\hs} \\
{\hb} &  & {\hb} & & {\hs} & \\
& {\hb} & & {\hb} & & {\hs} & & {\hs} \\
& & {\hb} & & {\hs} & & {\hs} & \\
& {\hs} & & {\hb}^y & & {\hs} & & {\hs} \\
} 
\]
By Lemma~\ref{l:mu.winverse}, we may assume that the entry $y'$ is not in the
heap by taking a distinct pair of $\mu$-masks for $x^{-1}$ on $w^{-1}$ if
necessary, because the heap of $w^{-1}$ is obtained from the heap of $w$ by
flipping the diagram upside down.  

Suppose there are three entries in the maximal northwest diagonal of the heap of
$w$ and consider the decorated heap associated to a $\mu$-mask on $w$.  If any
of the entries along the maximal northwest diagonal of the heap of $w$ are a
plain-zero then the entry must be the left critical zero of a zero-defect by
Lemma~\ref{l:pc}.  Hence, the decorated heap of any $\mu$-mask on $w$ has one of
the forms shown in Figure~\ref{f:mu_a1} on the entries in the maximal northwest
diagonal.

\begin{figure}[h]
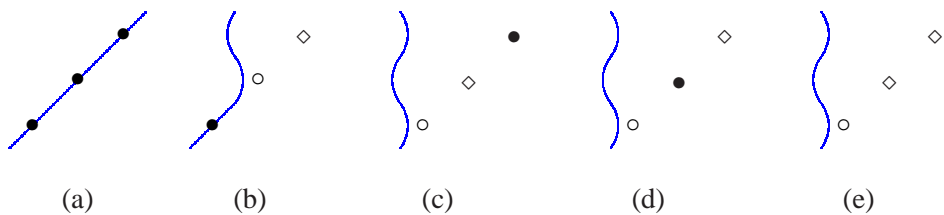

\begin{tabular}{ccccc}
\heap {
{\hs} & & \StringRX{\hf} \\
& \StringRX{\hf} \\
\StringRX{\hf} \\
 & & \\
} &
\heap {
\StringR{\hs} & & {\hd} \\
& \StringL{\hz} \\
\StringRX{\hf} \\
 & & \\
} &
\heap {
& \StringL{\hs} & & {\hf} \\
\StringR{\hs} & & {\hd} \\
& \StringL{\hz} \\
 & & \\
} &
\heap {
& \StringL{\hs} & & {\hd} \\
\StringR{\hs} & & {\hf} \\
& \StringL{\hz} \\
 & & \\
} &
\heap {
& \StringL{\hs} & & {\hd} \\
\StringR{\hs} & & {\hd} \\
& \StringL{\hz} \\
 & & \\
} \\
(a) & (b) & (c) & (d) & (e) \\
\end{tabular}
\caption{Maximal diagonal of Deodhar $\mu$-masks}
\label{f:mu_a1}
\end{figure}

In the arguments below, we will frequently use the fact that the strings on a
decorated heap associated to the mask $\s$ encode the 1-line notation for the
element $x = w^{\s}$.  If we label the strings according to their position along
the bottom of the decorated heap then the labeled strings must have the same
configuration of positions at the top of the decorated heap of $\s$ as they do
in the decorated heap of $\g$ because $w^{\s} = x = w^{\g}$.

{\bf Claim (1):  The $\mu$-masks $\s$ and $\g$ must fall into the same case (a) -
(e) of Figure~\ref{f:mu_a1}.}

Suppose $\s$ falls into case (a) and add strings to the decorated heap
associated to $\s$.  Since there is a single entry in the first column, the
string which is leftmost on the bottom of the decorated heap of $\s$ will end up
in the fourth position at the top of the decorated heap of $\s$.  Because
$w^{\s} = w^{\g}$ and none of the other cases send the string which is leftmost
on the bottom of the decorated heap to the fourth position on the top, we must
have that $\s$ falls into case (a) if and only if $\g$ falls into case (a).
Similarly, if $\s$ is in case (b) then the string which is leftmost on the
bottom of the decorated heap of $\s$ will end up in the second position at the
top of the decorated heap of $\s$, and we observe that none of the other cases
have this feature.  Hence, $\s$ falls into case (b) if and only if $\g$ falls
into case (b).

\begin{figure}[h]
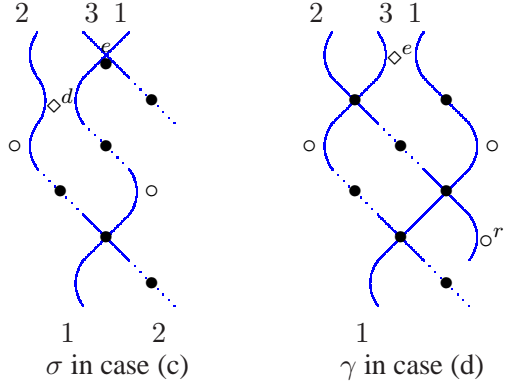

\begin{tabular}{cc}
\heap {
\ \ 2 & & 3 \ \ 1 & & \\
\StringR{\hs} & & \StringLRX{\stackrel{e}{\hf}} & & \\
& \StringLR{\hd^d} & & \StringLXD{\hf} & & \\
\StringR{\hz} & & \StringLXD{\hf} \\
& \StringLXD{\hf} & & \StringL{\hz} & & \\
{\hs} & & \StringLRX{\hf} \\
& \StringR{\hs} & & \StringLXD{\hf} & & \\
&   \ \ 1       & &  \ \ 2 \\
} &
\heap {
\ \ 2 & & 3 \ \ 1 & & \\
\StringR{\hs} & & \StringLR{\hd^e} & & \\
& \StringLRX{\hf} & & \StringLX{\hf} & & \\
\StringR{\hz} & & \StringLXD{\hf} & & \StringL{\hz} \\
& \StringLXD{\hf} & & \StringLRX{\hf} & & \\
{\hs} & & \StringLRX{\hf} & & \StringL{\hz^r} \\
& \StringR{\hs} & & \StringLXD{\hf} & & \\
&   \ \ 1       & &  \ \  \\
} \\
$\s$ in case (c) & $\g$ in case (d) \\
\end{tabular}
\caption{Case (c) distinguished from case (d)}
\label{f:case_cd}
\end{figure}

Next, suppose $\s$ falls into case (c) and $\g$ falls into case (d).  Then we
have decorated heap fragments of the form shown in Figure~\ref{f:case_cd}.
First, observe that the string labeled 2 in $\s$ cannot encounter a mask-value 0
entry along the minimal diagonal until after it has crossed the string labeled
1, since otherwise $d$ is not a zero-defect.  If the string labeled 3 in $\s$
crosses the string labeled 1 below $e$ then $e$ is a maximal one-defect in $\s$,
contradicting Lemma~\ref{l:no_isolation}.  Hence, the string labeled 3 must end
up lying to the right of the string labeled 1 on the bottom of the decorated
heap of $\g$ since $w^{\s} = w^{\g}$.  This implies that the string labeled 3
cannot encounter a mask-value 0 entry along the minimal diagonal in $\g$ until
after it has crossed the string labeled 1.  Since the string labeled 2 lies to
the right of the string labeled 1 on the bottom of the decorated heap of $\s$,
we have that the string labeled 2 cannot encounter a mask-value 0 entry from the
top in $\g$ until after it has crossed the string labeled 1.  Therefore, the
position of the right critical zero of $e$ in $\g$ is completely determined by
the position of the string labeled 1 on the bottom of the decorated heap of
$\s$.  Furthermore, the string labeled 2 must eventually encounter a mask-value
0 entry $r$ in $\g$ so that it lies in the same position on the bottom of the
decorated heap of $\g$ as it does in $\s$, since $w^{\s} = w^{\g}$.  Indeed,
lateral convexity implies that if the string labeled 2 does not encounter a
mask-value 0 entry in $\g$ then the string ends up on the bottom in a column
that is strictly right of any column that can be reached by a string traveling
along the minimal southwest diagonal of the heap of $w$.

By Lemma~\ref{l:pc}, $r$ is either a left critical zero of a zero-defect $f$ or
$r$ is a zero-defect in $\g$, and so there exists an entry directly southeast
of $r$ to facilitate the string crossing for the zero-defect.  Thus, there can
be no entry northeast of the right critical zero of $e$, for otherwise we
obtain an I-shape in the heap of $w$, contradicting that $w$ is Deodhar.
However, since the right critical zero of $e$ is not a zero-defect in $\g$
because the path of the string labeled 1 is prescribed, this implies that $e$
and one of $\{f,r\}$ form a separating pair for $\g$, contradicting that $\g$
is a $\mu$-mask by Lemma~\ref{l:connected}.

Next, suppose $\s$ falls into case (e) as illustrated in Figure~\ref{f:case_e}.
Observe that the 1-line notation for $x = w^{\s}$ determines the location of the
right critical zeros of the zero-defects of $\s$.
\begin{figure}[h]
\begin{tabular}{c}
\heap {
 \ \ 3 & & 2 \ \ 1 & & \\
\StringR{\hs} & & \StringLR{\hd^e} & & \\
& \StringLR{\hd^d} & & \StringL{\hz^p} & & \\
\StringR{\hz} & & \StringLXDR{\hf} & & \\
& \StringLXDR{\hf} & & \StringL{\hz^q} & & \\
\StringR{\hs} & & \StringLRX{\hf} & & \\
 \ \ 1 & & 2 \ \ 3 \\
}
\end{tabular}
\caption{$\s$ in case (e)}
\label{f:case_e}
\end{figure}
To see this, first observe that the left string of $e$ becomes the right string
of $d$ in $\s$, and the strings of $e$ must cross above the right critical zero
$q$ of $d$.  Otherwise we obtain a 4-stack, contradicting that $w$ is Deodhar.
Hence, all of the entries along the path of the right string of $e$ have
mask-value 1 below the right critical zero $p$ of $e$.  Thus, the position of
the string labeled 1 on the bottom of the heap determines the column where $p$
appears.  Similarly, the position of the string labeled 2 on the bottom of the
heap determines the column where $q$ appears.  Since $w^{\s} = w^{\g}$, these
strings must occur in the same positions at the top and bottom of the decorated
heap of $\g$.

\begin{figure}[h]
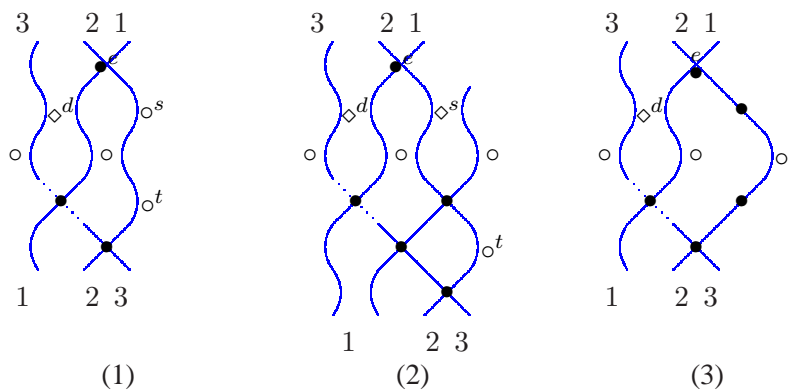

\begin{tabular}{ccc}
\heap {
 \ \ 3 & & 2 \ \ 1 & & \\
\StringR{\hs} & & \StringLRX{\hf^e} & & \\
& \StringLR{\hd^d} & & \StringL{\hz^s} & & \\
\StringR{\hz} & & \StringLR{\hz} & & \\
& \StringLXDR{\hf} & & \StringL{\hz^t} & & \\
\StringR{\hs} & & \StringLRX{\hf} & & \\
 \ \ 1 & & 2 \ \ 3 \\
} &
\heap {
 \ \ 3 & & 2 \ \ 1 & & \\
\StringR{\hs} & & \StringLRX{\hf^e} & & \\
& \StringLR{\hd^d} & & \StringLR{\hd^s} & & \\
\StringR{\hz} & & \StringLR{\hz} & & \StringL{\hz} \\
& \StringLXDR{\hf} & & \StringLRX{\hf} & & {\hs} \\
\StringR{\hs} & & \StringLRX{\hf} & & \StringL{\hz^t} \\
& \StringLR{\hs} & & \StringLRX{\hf} & & {\hs} \\
& 1 \ \ & & 2 \ \ 3 \\
} &
\heap {
 \ \ 3 & & 2 \ \ 1 & & \\
\StringR{\hs} & & \StringLRX{\stackrel{e}{\hf}} & & \\
& \StringLR{\hd^d} & & \StringLX{\hf} & & \\
\StringR{\hz} & & \StringL{\hz} & & \StringL{\hz^s}\\
& \StringLXDR{\hf} & & \StringRX{\hf} & & \\
\StringR{\hs} & & \StringLRX{\hf} & & \\
 \ \ 1 & & 2 \ \ 3 \\
} \\
(1) & (2) & (3) \\
\end{tabular}
\caption{$\g$ in case (c) distinguished from $\s$ in case (e)}
\label{f:case_ce}
\end{figure}

If $\g$ falls into case (c) then the decorated heap of $\g$ has one of the forms
shown in Figure~\ref{f:case_ce}.  To see this, begin by observing that the
string labeled 3 in $\g$ must cross both of the strings labeled 1 and 2 since
$w^{\s} = w^{\g}$.  Therefore, the string labeled 3 does not encounter a
mask-value 0 entry along the minimal diagonal until it crosses both the strings
labeled 1 and 2.  Hence, the location of the right critical zero of $d$ in $\g$
is completely determined by the location of the string labeled 1 on the bottom
of the decorated heap of $\s$.  Also, observe that the string labeled 2 in $\g$
cannot encounter a mask-value 0 entry until the column after the right critical
zero of $d$ appears, for otherwise the string labeled 2 ends up lying to the
left of the string labeled 1 on the bottom of the decorated heap of $\g$,
contradicting that $w^{\s} = w^{\g}$.  Since the string labeled 2 must cross the
string labeled 3 in $\g$, the string labeled 2 eventually encounters a
mask-value 0 entry $s$, and the string labeled 2 either turns once or the string
labeled 2 turns three times before crossing the string labeled 3.

If the string labeled 2 turns three times in the heap before crossing the string
labeled 3, and $s$ is a plain-zero, then $\g$ has the form shown in
Figure~\ref{f:case_ce}(1).  
In this case, $t$ is a plain-zero so there exists a
zero-defect with $t$ as its left critical zero, and the entry directly southeast
of $t$ lies in the heap of $w$.  Also, $s$ must be the left critical zero of a
zero-defect by Lemma~\ref{l:pc}, which introduces an I-shape into the heap
contradicting that $w$ is Deodhar.

If the string labeled 2 turns three times in the heap before crossing the string
labeled 3, and $s$ is a zero-defect, then $\g$ has the form shown in
Figure~\ref{f:case_ce}(2).  Observe that if the right string of $s$ does not
cross the left string of $s$ before the third turning point $t$ of the string
labeled 2 then the strings of $s$ do not have an opportunity to meet again so
$s$ cannot be a zero-defect.  Moreover, $t$ is a plain-zero so there exists a
zero-defect $f$ with $t$ as its left critical zero, and the entry directly
southeast of $t$ lies in the heap of $w$.  Also, the right critical zero of $s$
is a plain-zero, and $s$ cannot share a critical zero with another zero-defect
to the right because this would introduce an I-shape, contradicting that $w$ is
Deodhar.  Hence, $s$ and $f$ form a separating pair for $\g$ contradicting that
$\g$ is a $\mu$-mask by Lemma~\ref{l:connected}.

If the string labeled 2 turns once in the heap before crossing the string
labeled 3 then $\g$ has the form shown in Figure~\ref{f:case_ce}(3).  In this
case, the right critical zero of $d$ is a plain-zero which is not shared with
another zero-defect to the right, because this would introduce an I-shape into
the heap of $w$, contradicting that $w$ is Deodhar.  However, $s$ is a
plain-zero which must be the left critical zero of a zero-defect $g$ by
Lemma~\ref{l:pc}, so $g$ and $d$ form a separating pair for $\g$,
contradicting that $\g$ is a $\mu$-mask by Lemma~\ref{l:connected}.

Thus, we have shown that $\g$ lies in case (c) if and only if $\s$ lies in case
(c).  Finally, if $\g$ falls into case (d) and $\s$ falls into case (e) then the
decorated heap of $\g$ has form shown in Figure~\ref{f:case_de}.

\begin{figure}[h]
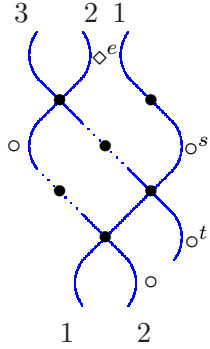

\begin{tabular}{c}
\heap {
 \ \ 3 & & 2 \ \ 1 & & \\
\StringR{\hs} & & \StringLR{\hd^e} & & \\
& \StringLRX{\hf} & & \StringLX{\hf} & & \\
\StringR{\hz} & & \StringLXD{\hf} & & \StringL{\hz^s}\\
& \StringLXD{\hf} & & \StringLRX{\hf} & & \\
{\hs} & & \StringLRX{\hf} & & \StringL{\hz^t} \\
& \StringR{\hs} & & \StringL{\hz} & & \\
&  \ \ 1        & & 2 \ \ \\
}
\end{tabular}
\caption{$\g$ in case (d) distinguished from $\s$ in case (e)}
\label{f:case_de}
\end{figure}

Since $w^{\s} = w^{\g}$, the string labeled 2 in $\g$ cannot encounter a
mask-value 0 entry along the minimal diagonal until it has crossed the string
labeled 1.  Hence, the string labeled 3 cannot encounter a mask-value 0 entry
from the top until it has crossed the string labeled 1.  Therefore, the position
of the right critical zero $s$ of $e$ is completely determined by the position
of the string labeled 1 on the bottom of the decorated heap of $\s$.  Since the
string labeled 3 in $\s$ travels along the minimal southwest diagonal in $\s$,
we have that the string labeled 3 in $\g$ must encounter a mask-value 0 entry
$t$ in $\g$ sometime after it has crossed the string labeled 1.  Indeed, lateral
convexity implies that if the string labeled 3 does not encounter a mask-value 0
entry in $\g$ then the string ends up on the bottom in a column that is strictly
right of any column that can be reached by a string traveling along the minimal
southwest diagonal of the heap of $w$.  The entry $t$ is either a zero-defect or
the left critical zero of a zero-defect $h$ by Lemma~\ref{l:pc}, so the entry
directly southeast of $t$ lies in the heap of $w$ to facilitate the string
crossing for a zero-defect.  Hence, there is no entry northeast of $s$, for
otherwise we introduce an I-shape into the heap of $w$, contradicting that $w$
is Deodhar.  Since the right critical zero of $e$ is a plain-zero, this implies
that $e$ does not share a critical zero with any zero-defect lying to the right
of $e$, so $e$ and one of $\{t, h\}$ form a separating pair for $\g$,
contradicting that $\g$ is a $\mu$-mask by Lemma~\ref{l:connected}.

Thus, we have shown that $\s$ and $\g$ fall into the same case from
Figure~\ref{f:mu_a1}, completing the proof of Claim (1).

{\bf Claim (2):  For each case (a) - (e) of Figure~\ref{f:mu_a1}, we may
construct a smaller length element with distinct $\mu$-masks, contradicting our
minimal length choice of $w$.}

In cases (a), (b), (c) and (d), we can simply remove the three entries on the
maximal northwest diagonal from the decorated heaps of $\s$ and $\g$.  Since
$\g$ and $\s$ fall into the same case, the mask-values agree on these entries.
Because the entries are maximal in the decorated heap, removing them does not
alter the defect status of any other entries in the decorated heap.  Also, in
each case we remove the same number of zero-defects as plain-zeros.  Hence, we
obtain distinct $\mu$-masks $\check{\s}$ and $\check{\g}$ on a smaller length
element $\check{w}$ because $\s$ and $\g$ were assumed to be distinct.  Finally,
since $w^{\s} = w^{\g}$ and the masks agree on the maximal entries that we
removed, we have that $\check{w}^{\check{\s}} = \check{w}^{\check{\g}}$.

Consider the case that $\s$ and $\g$ are in case (e).  Recall from the
discussion of Figure~\ref{f:case_e} that the strings of the zero-defect $e$
always cross above the right critical zero $q$ for the zero-defect $d$, and
1-line notation of $x = w^{\s} = w^{\g}$ completely determines the location of
the right critical zeros $p$ and $q$ of the zero-defects $e$ and $d$.  In
particular, $p$ and $q$ each occur at the same entry of the heap of $w$ in both
$\s$ and $\g$.  If entries exist in the heap of $w$ directly northeast of $p$
and directly southeast of $q$ then we obtain an I-shape in the heap of $w$,
contradicting that $w$ is Deodhar.

If $p$ has no entries lying to the northeast then the same strings encounter
the entry $p$ from the top in the decorated heaps of both $\s$ and $\g$.  Hence,
we may remove $e$ and change the mask-value of $p$ to 1 in $\s$ and $\g$.  This
removes one zero-defect and one plain-zero from the masks, so we obtain
$\mu$-masks $\check{\s}$ and $\check{\g}$ on a smaller length element
$\check{w}$, and the masks are distinct because $\s$ and $\g$ were assumed to be
distinct.  Since $w^{\s} = w^{\g}$ and the same strings encounter $p$ in both
masks, we have that $\check{w}^{\check{\s}} = \check{w}^{\check{\g}}$.

If $q$ has no entries lying to the southeast then the same strings encounter
the entry $q$ from the bottom in the decorated heaps of both $\s$ and $\g$.
Also, $d$ does not share a critical-zero with any zero-defect lying to the
right.  Hence, we may change the mask-value of the left critical zero of $d$ to
1, and change the mask-value of the right critical zero $q$ of $d$ to 1 in $\s$
and $\g$.  This turns $d$ into a plain-zero, but does not change the defect
status of any of the other entries in the decorated heaps.  Thus, we obtain
masks $\check{\s}$ and $\check{\g}$ on $w$ where $\check{\s}$ and $\check{\g}$
fall into case (b) of Figure~\ref{f:mu_a1}.  Also, $D(\s) = D(\check{\s}) = 1 =
D(\check{\g}) = D(\g)$ so $\check{\s}$ and $\check{\g}$ are $\mu$-masks.  The
masks are distinct because $\s$ and $\g$ were assumed to be distinct.  Since
$w^{\s} = w^{\g}$ and the same strings encounter the entries that we changed in
both masks, we have that $\check{w}^{\check{\s}} = \check{w}^{\check{\g}}$.  By
applying the argument for case (b) above, we can then remove the maximal
northwest diagonal of the decorated heaps associated to $\check{\s}$ and
$\check{\g}$ to obtain a smaller length element with distinct $\mu$-masks.

Thus, in all cases we contradict our minimal length choice of $w$.  This
completes the proof of Claim (2).

Finally, suppose there are less than three entries in the maximal northwest
diagonal of the heap of $w$.  If there is a single entry in the maximal
northwest diagonal then any $\mu$-mask must have mask-value 1 on this entry by
Lemma~\ref{l:pc}.  If there are two entries in the maximal northwest diagonal of
the decorated heap of $w$ associated to $\s$ then either both entries have
mask-value 1, or the maximal entry of the diagonal is a zero-defect and the
minimal entry is a plain-zero.  Since the leftmost string in the decorated heap
corresponds to a fixed point in the permutation $x = w^{\s}$ if and only if the
mask-value of the leftmost entry is 0, we see that the mask-values of the
entries in the maximal northwest diagonal agree in the decorated heaps
associated to $\s$ and $\g$.  Therefore, we can remove the entries along the
maximal northwest diagonal in both decorated heaps.  Doing so removes the same
number of plain-zeros as zero-defects so the resulting masks are distinct
$\mu$-masks on a smaller length element that encode the same element,
contradicting our minimal choice of $w$.
\end{proof}

\begin{corollary}\label{c:main.b}
If $w$ is a Deodhar element of type $B$ then $\mu(x,w) \in \{0,1\}$ for all $x$.
\end{corollary}
\begin{proof}
Let $w$ be a Deodhar element of type $B$.  By \cite[Corollary 5.3]{b-j1}, the
Deodhar elements in type $B$ are short-braid avoiding and an analogue of the
Lateral Convexity lemma \cite[Lemma 1]{b-w} for type $B$ shows that there can be
at most one $s_0$ generator in $w$.  Therefore, we can interpret any reduced
expression for $w$ as a type $A$ reduced expression by sending the Coxeter
generators $s_0, s_1, \ldots, s_{n-1}$ of type $B$ to the Coxeter generators
$s_1, s_2, \ldots, s_n$ of type $A$, even though $(s_0 s_1)^3 \neq 1$ in type
$B$.  This enables us to view the decorated heap of a Deodhar type $B$ element
as the decorated heap of a Deodhar permutation.  Moreover, it follows quickly
from the discussion in Section~\ref{s:type_bd} that an entry $z$ is a
zero-defect in the decorated heap of type $B$ if and only if $z$ is a
zero-defect in the corresponding decorated heap of type $A$.  In fact, this
shows that the Deodhar decorated heaps of type $B_n$ are in bijective
correspondence with the Deodhar decorated heaps of $A_n$ and the bijection
preserves the Deodhar statistic.  Hence, if there exist distinct $\mu$-masks
$\s$ and $\g$ on $w$ such that $w^{\s} = w^{\g}$ then the bijection gives a
contradiction to Theorem~\ref{t:mu.main}.
\end{proof}


\bigskip
\section{The 0-1 result for type $D$ Deodhar elements}\label{s:mu.d}

In this section we sketch the proof of Theorem~\ref{t:mu.remain} for the Deodhar
elements of type $D$ using the methods developed in \cite{b-j1} to characterize
the Deodhar elements of type $D$.  In particular, we treat the convex and
non-convex elements separately.  The proof for the convex elements closely
follows the proof of our classification result \cite[Proposition 8.7]{b-j1} and
we cite relevant facts from that proof as necessary.

We begin by recalling a map which is used to obtain a reduction to type $A$.
\begin{lemma}\label{l:d.pi}
Let $w$ be a convex Deodhar element of type $D$ and suppose $\s$ is a mask on
$w$ with three entries in column 1 and no entry directly northeast of the top
entry in column 1.  Then, there exists a map $\pi$ which we denote $(w, \s)
\mapsto (\pi(w), \pi(\s))$ with the following properties:
\begin{enumerate}
\item[(1)]  $\pi(w)$ is a Deodhar element of type $A$.
\item[(2)]  $\pi(\s)$ is a mask on $\pi(w)$.
\item[(3)]  $D(\pi(\s)) = D(\s)$.
\item[(4)]  The entries lying to the right of column 3 in the heap of $\pi(w)$
correspond to the entries lying to the right of column 1 in the heap of $w$.
Moreover, the mask-value and defect-status of these entries in decorated heap of
$\pi(\s)$ agree with the mask-value and defect-status of the corresponding
entries in the decorated heap of $\s$.
\end{enumerate}
\end{lemma}
\begin{proof}
The map of elements is defined in \cite[Definition 8.3]{b-j1}.  This definition
extends to a map of masks as described in Appendix~\ref{s:appendix} with the
given properties by considering the proof of \cite[Proposition 8.7]{b-j1}.
\end{proof}

\begin{proposition}\label{p:mu.d.convex}
Suppose $w$ is a convex Deodhar type $D$ element, and $x$ is any element of type
$D$.  Then, $\mu(x, w) \in \{0, 1\}$.
\end{proposition}

\begin{proof}
Suppose for the sake of contradiction that there exist distinct $\mu$-masks $\s$
and $\g$ such that $w^{\s} = x = w^{\g}$.  Let $k(w)$ be the number of distinct
levels in the first column of the coalesced heap of $w$.  We consider the same
cases that were considered in the proof of \cite[Proposition 8.7]{b-j1}.  

\bigskip

\noindent \textbf{Case $\mathbf{k(w)=0}$ or $\mathbf{1}$}.  
If $w$ has at most a single entry in the first column then $w$ is contained in a
parabolic subgroup of type $A$.  Hence, Theorem~\ref{t:mu.main} gives that
$\mu(x,w) \in \{0, 1\}$ for all $x$.

Suppose there exist two entries in column 1.  If exactly one of the entries in
column 1 has mask-value 0 in $\s$ or $\g$, then the mask-value 0 entry cannot be
the left critical zero for any zero-defect because the left string of such a
zero-defect will be negative while the right string will be positive, and this
would contradict Lemma~\ref{l:no_isolation}.  Hence, both entries in column 1
must have the same mask-value in any $\mu$-mask.  Moreover, both entries of
column 1 have mask-value 0 in $\s$ if and only if all of the entries in the
1-line notation for $x = w^{\s}$ are positive, and this occurs if and only if
both entries of column 1 have mask-value 0 in $\g$ because $w^{\s} = x =
w^{\g}$.

Next, consider the decorated heaps $\check{\s}$ and $\check{\g}$ obtained from
$\s$ and $\g$ by removing the plain-zero $s_{\tilde{1}}$ entry from column 1 of
each heap.  Since $\s$ and $\g$ were assumed to be $\mu$-masks in type $D$ for a
convex element, we have that $\check{\s}$ and $\check{\g}$ are masks on a
short-braid avoiding permutation $\check{w}$.

If $\s$ and $\g$ have mask-value 1 for both entries in column 1 and $\check{w}$
is Deodhar then $\check{\s}$ and $\check{\g}$ contradict
Theorem~\ref{t:mu.main}.  Also, if $\s$ and $\g$ have mask-value 0 for both
entries in column 1 then the masks $\check{\s}$ and $\check{\g}$ have Deodhar
statistic 0.  Hence, we may assume $\check{w}$ is a non-Deodhar permutation.  By
\cite[Theorem 1]{b-w} this implies that the heap of $\check{w}$ contains a
hexagon that uses the single entry in column 1.  By \cite[Theorem 8.1]{b-j1},
the only convex Deodhar elements in type $D$ which project in this way to type
$A$ are:
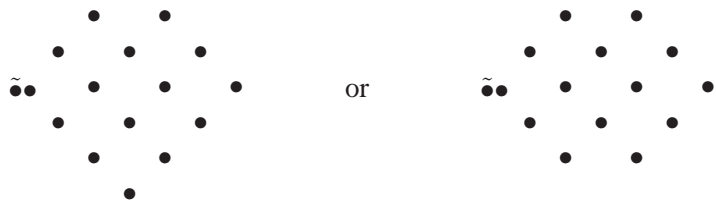
\[
\xymatrix @=-4pt @! {
& {\hs} & & {\hs} & & {\hs} & & {\hs} \\
& & {\hf} & & {\hf} & & {\hs} & \\
& {\hf} & & {\hf} & & {\hf} & & {\hs} \\
\tilde{\hf} {\hf} &  & {\hf} & & {\hf} & & {\hf} \\
& {\hf} & & {\hf} & & {\hf} & & {\hs} \\
& & {\hf} & & {\hf} & & {\hs} & \\
& {\hs} & & {\hf} & & {\hs} & & {\hs} \\
} 
\parbox[t]{1in}{ \vspace{0.5in} \hspace{.3in} or }
\xymatrix @=-4pt @! {
& {\hs} & & {\hs} & & {\hs} & & {\hs} \\
& & {\hf} & & {\hf} & & {\hs} & \\
& {\hf} & & {\hf} & & {\hf} & & {\hs} \\
\tilde{\hf} {\hf} &  & {\hf} & & {\hf} & & {\hf} \\
& {\hf} & & {\hf} & & {\hf} & & {\hs} \\
& & {\hf} & & {\hf} & & {\hs} & \\
& {\hs} & & {\hs} & & {\hs} & & {\hs} \\
} 
\]
It is straightforward to verify that $\mu(x, w) \in \{0,1\}$ for these elements.
There is C++ code available at
\url{http://www.math.ucdavis.edu/~brant/code/} for this purpose.

\bigskip

\noindent \textbf{Case $\mathbf{k(w)=2}$}.
By \cite[Lemma 6.7]{b-j1}, whenever the first column of a type $D$ heap contains
entries on more than one distinct level, we must have that each level of column
1 contains a unique entry, and the entries alternate between the $s_1$ and
$s_{\tilde{1}}$ generators.  Without loss of generality, we can assume assume
the top entry in column 1 of the heap of $w$ is $s_{\tilde{1}}$ since $\mu(x,w)$
is invariant under applying a Coxeter graph automorphism to $x$ and $w$.  Since
$w$ is Deodhar, we have by Lemma~\ref{l:shape} that there are at most 3 entries
in column 3.  By considering $w^{-1}$ if necessary and applying
Lemma~\ref{l:mu.winverse}, we can further assume that the maximal northwest
diagonal from the top entry in column 1 contains at most two entries:
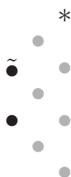
\[
\xymatrix @=-4pt @! {
& {\hs} & & {\hv} & & {\hs} & & {\hs} \\
& & {\hb} & & {\hs} & & {\hs} & \\
& {\tilde{\hf}} & & {\hb} & & {\hs} & & {\hs} \\
{\hs} &  & {\hb} & & {\hs} & \\
& {\hf} & & {\hb} & & {\hs} & & {\hs} \\
& & {\hb} & & {\hs} & & {\hs} & \\
& {\hs} & & {\hb} & & {\hs} & & {\hs} \\
} 
\]

Consider the decorated heap associated to a $\mu$-mask on $w$.  The maximal
northwest diagonal of any $\mu$-mask on $w$ has one of the forms shown in
Figure~\ref{f:mu_k2} by Lemma~\ref{l:no_isolation}.

\begin{figure}[h]
\begin{tabular}{ll}
\xymatrix @=-4pt @! {
& {\hd} & & {\hs} & & \\
\tilde{\hz} & & \\
 & & \\
} &
\xymatrix @=-4pt @! {
& {\hf} & & {\hs} & & \\
\tilde{\hf} & & \\
 & & \\
} \\
(a) & (b) \\
\end{tabular}
\caption{$\mu$-masks with two levels in column 1}
\label{f:mu_k2}
\end{figure}
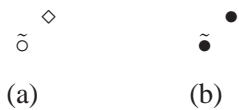

Observe that some strings are labeled negatively in case (b) but all strings are
labeled positively in case (a).  Since $w^{\s} = x = w^{\g}$, we have that both
masks must fall into the same case, so the mask-values of $\s$ and $\g$ agree
along the maximal northwest diagonal.  As in the proof of \cite[Proposition
8.7]{b-j1}, we may remove this northwest diagonal from $\s$ and $\g$ to obtain a
pair of decorated heaps $\check{\s}$ and $\check{\g}$ on a Deodhar type $A$
element $\check{w}$.  Observe that $\check{w}^{\check{\s}} =
\check{w}^{\check{\g}}$ and $D(\check{\s}) = D(\check{\g})$.  Since $\s$ and
$\g$ are distinct, we obtain a contradiction to Theorem~\ref{t:mu.main}.

\bigskip

\noindent \textbf{Case $\mathbf{k(w)=3}$}.  Suppose that $w$ has entries on 3
distinct levels in the first column and assume without loss of generality that
the top entry corresponds to $s_{1}$.  As in the proof of \cite[Proposition
8.7]{b-j1}, our strategy will be to project the heap of $w$ to a Deodhar type
$A$ element by adding the three entries marked as $\hh$ to the left of the heap
of $w$:
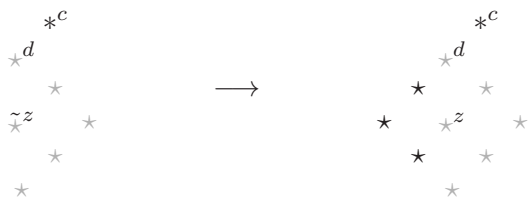
\[
\xymatrix @=-4pt @! {
& {\hs} & & {\hv^c} & & {\hs} & & {\hs} \\
& & {\ha}^d & & {\hs} & & {\hs} & \\
& {\hs} & & {\ha} & & {\hs} & & {\hs} \\
{\hs} & & \tilde{\ha}^z &  & {\ha} & & {\hs} & \\
& {\hs} & & {\ha} & & {\hs} & & {\hs} \\
& & {\ha} & & {\hs} & & {\hs} & \\
& {\hs} & & {\hs} & & {\hs} & & {\hs} \\
}
\begin{array}{c}
\\ \\ \\ \\
\longrightarrow
\end{array}
\hspace{.5in}
\xymatrix @=-4pt @! {
& {\hs} & & {\hv^c} & & {\hs} & & {\hs} \\
& & {\ha}^d & & {\hs} & & {\hs} & \\
& {\hh} & & {\ha} & & {\hs} & & {\hs} \\
{\hh} & & {\ha}^z &  & {\ha} & & {\hs} & \\
& {\hh} & & {\ha} & & {\hs} & & {\hs} \\
& & {\ha} & & {\hs} & & {\hs} & \\
& {\hs} & & {\hs} & & {\hs} & & {\hs} \\
}
\]
If $w$ contains a 4-stack in column 2 then $w$ is not Deodhar by
Lemma~\ref{l:shape}.  Hence, we can assume there are at most 3 entries in column
2.  By considering $w^{-1}$ if necessary and applying Lemma~\ref{l:mu.winverse},
we can assume that the point marked $c$ in column 2 is not in the heap.

Consider the decorated heap associated to a $\mu$-mask on $w$.  By
Lemma~\ref{l:no_isolation} every plain-zero must be associated with a
zero-defect and there are no maximal one-defects.  Hence, any $\mu$-mask on $w$
has one of the following forms:
\begin{enumerate}
\item[(a)]  $d$ is a plain-one.
\item[(b)]  $d$ is a zero-defect and $z$ has mask-value 1.
\item[(c)]  $d$ is a zero-defect and $z$ has mask-value 0.  The strings for $d$
meet at $z$.
\item[(d)]  $d$ is a zero-defect and $z$ has mask-value 0.  The strings for $d$
do not meet at $z$.
\end{enumerate}
Our first goal is to show that the $\mu$-masks $\s$ and $\g$ must fall into the
same case of (a) - (d).

Begin by observing that we cannot have $\s$ in case (b), and $\g$ in case (c) or
(d).  This follows because $w^{\s} = x = w^{\g}$ and $x$ has some negatively
signed entries in case (b), but $x$ has all positively signed entries in cases
(c) and (d).

Next, suppose that $\s$ is in case (b) and $\g$ is in case (a).  Then, we have
the decorated heap fragments shown in Figure~\ref{f:mu_k3ab}.

\begin{figure}[h]
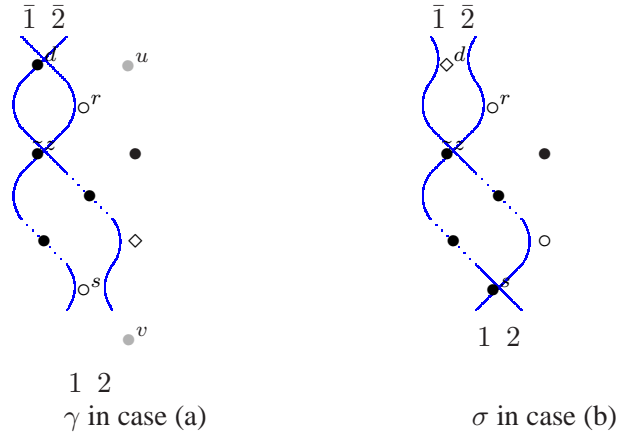

\begin{tabular}{cc}
\heap {
& {\hs} & \bar{1} \ \ \bar{2} & {\hs} & & {\hs} & & {\hs} \\
& & \StringLRX{\hf^d} & & {\hb^u} & & {\hs} & \\
& \StringR{\hs} & & \StringL{\hz^r} & & {\hs} & & {\hs} \\
{\hs} & & \StringLRX{\tilde{\hf}^z} &  & {\hf} & & {\hs} & \\
& \StringR{\hs} & & \StringLXD{\hf} & & {\hs} & & {\hs} \\
& & \StringLXD{\hf} & & \StringL{\hd} & & {\hs} & \\
& {\hs} & & \StringLR{\hz^s} & & {\hs} & & {\hs} \\
& & {\hs} & & {\hb^v} & & {\hs} & \\
& & {\hs} & 1 \ \ 2 & {\hs} & & {\hs} & & {\hs} \\
} &
\heap {
& {\hs} & \bar{1} \ \ \bar{2} & {\hs} & & {\hs} & & {\hs} \\
& & \StringLR{\hd^d} & & {\hs} & & {\hs} & \\
& \StringR{\hs} & & \StringL{\hz^r} & & {\hs} & & {\hs} \\
{\hs} & & \StringLRX{\tilde{\hf}^z} &  & {\hf} & & {\hs} & \\
& \StringR{\hs} & & \StringLXD{\hf} & & {\hs} & & {\hs} \\
& & \StringLXD{\hf} & & \StringL{\hz} & & {\hs} & \\
& {\hs} & & \StringLRX{\hf^s} & & {\hs} & & {\hs} \\
& & {\hs} & 1 \ \ 2 & {\hs} & & {\hs} & & {\hs} \\
} \\
$\g$ in case (a) & $\s$ in case (b) \\
\end{tabular}
\caption{$\mu$-masks with three levels in column 1}
\label{f:mu_k3ab}
\end{figure}

The justification for Figure~\ref{f:mu_k3ab} runs as follows.  In case (b), the
left string of $d$ must become negative as it passes through the $s_{\tilde{1}}$
entry.  If we label the string of $d$ which is leftmost on the bottom of the
decorated heap by 1, and we label the string of $d$ which is rightmost on the
bottom of the decorated heap by 2, then in order for $d$ to be a defect the
strings of $d$ must be labeled at the top by $\bar{1} \bar{2}$.  This forces the
right string of $d$ to also pass through the $s_{\tilde{1}}$ entry so the right
critical zero $r$ of $d$ lies next to $d$ as shown.  Moreover, the string
labeled 2 does not encounter a mask-value 0 entry below $z$ in $\s$ until it
crosses the string labeled 1.  Since $w^{\g} = w^{\s}$, we must have these
strings labeled in the same way in $\g$.  Hence, the entry $r$ has mask-value 0
in $\g$.  The entry $r$ is not a zero-defect because the left string of $r$ is
negatively labeled while the right string of $r$ is positively labeled.  Also,
the string labeled 1 must encounter a mask-value 0 entry at $s$ in $\g$ in order
for the strings of $d$ to be labeled in $\g$ as they are in $\s$.  The entry $s$
is a plain-zero in $\g$, for otherwise we introduce a 4-stack into the heap,
which contradicts that $w$ is Deodhar.  Therefore, there must be zero-defects
which have $r$ and $s$ as left critical zeros by Lemma~\ref{l:no_isolation}.
But this implies that $u$ and $v$ appear the decorated heap of $\g$, which
introduces an I-shape into the heap of $w$, contradicting that $w$ is Deodhar.
Thus, we have shown that $\s$ is in case (b) if and only if $\g$ is in case (b).

\begin{figure}[h]
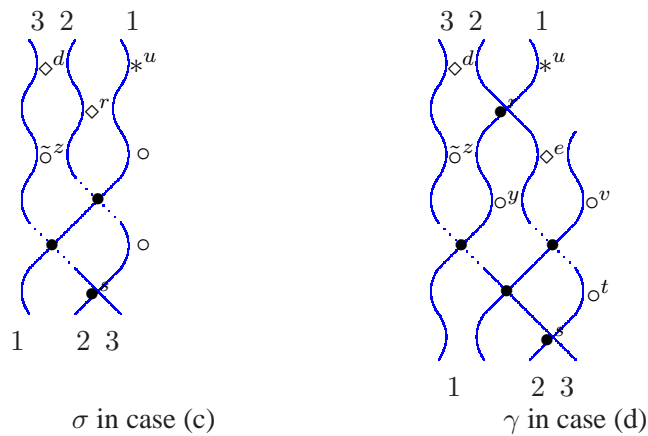

\begin{tabular}{cc}
\heap {
& {\hs} & 3 \ \ 2 & & 1 \ \ \  & {\hs} & & {\hs} \\
& & \StringLR{\hd^d} & & \StringL{\hv^u} & & {\hs} & \\
& \StringR{\hs} & & \StringLR{\hd^r} & & {\hs} & & {\hs} \\
{\hs} & & \StringLR{\tilde{\hz}^z} & & \StringL{\hz} & & {\hs} & \\
& \StringR{\hs} & & \StringLXDR{\hf} & & {\hs} & & {\hs} \\
& & \StringLXDR{\hf} & & \StringL{\hz} & & {\hs} & \\
& \StringR{\hs} & & \StringLRX{\hf^s} & & {\hs} & & {\hs} \\
& \ \ \ 1 & & 2 \ \ 3 & {\hs} & & {\hs} & & {\hs} \\
} &
\heap {
& {\hs} & 3 \ \ 2 & & 1 \ \ \ & {\hs} & & {\hs} \\
& & \StringLR{\hd^d} & & \StringL{\hv^u} & & {\hs} & \\
& \StringR{\hs} & & \StringLRX{\hf^r} & & {\hs} & & {\hs} \\
{\hs} & & \StringLR{\tilde{\hz}^z} &  & \StringLR{\hd^e} & & {\hs} & \\
& \StringR{\hs} & & \StringLR{\hz^y} & & \StringL{\hz^v} & & {\hs} \\
& & \StringLXDR{\hf} & & \StringLXDR{\hf} & & {\hs} & \\
& \StringR{\hs} & & \StringLRX{\hf} & & \StringL{\hz^t} & & {\hs} \\
& & \StringLR{\hs} & & \StringLRX{\hf^s} & & {\hs} & \\
& {\hs} & 1 \ \  & {\hs} & 2 \ \ 3 & &  & {\hs} & & {\hs} \\
} \\
$\s$ in case (c) & $\g$ in case (d)
\end{tabular}
\caption{$\mu$-masks with three levels in column 1}
\label{f:mu_k3cd}
\end{figure}

Next, suppose that $\s$ is in case (c) and $\g$ is in case (d).  Then, we have
the decorated heap fragments shown in Figure~\ref{f:mu_k3cd}.  In case (c),
observe that the right critical zero $r$ for $d$ must lie next to $d$ as shown
in order for the strings of $d$ to meet at $z$.  Also, the string labeled 3
cannot encounter a mask-value 0 entry below $z$ in $\s$ until it crosses the
string labeled 2 at some point $s$ on the minimal southwest diagonal because $d$
is a defect.  Note that $z$ is a plain-zero as there is no other $s_{\tilde{1}}$
generator in the heap of $w$ to support a defect at $z$, and if we change the
mask-value of $z$ to 1 then $d$ remains a zero-defect.  Hence, $r$ must be a
zero-defect with $z$ as its left critical zero by Lemma~\ref{l:no_isolation}.
The right string of $d$ becomes the left string of $r$ and the strings of $r$
must cross above the point directly northeast of $s$.  Otherwise the strings of
$r$ will not have an opportunity to cross again.  Hence, $s$ occurs in column 2
or further to the right, and so $u$ is not in the heap as it introduces an
I-shape.

Observe that the string which is in the third position at the top of $\s$ ends
up on the bottom of $\s$ in a position lying to the left of the strings labeled
2 and 3.  Since $w^{\s} = w^{\g}$ the strings in $\g$ must have the same
positions at the top and bottom of the decorated heap as in $\s$.  The string
labeled 3 cannot encounter a mask-value 0 entry below $z$ in $\g$ until it
crosses the string labeled 2.  Since the strings of $d$ do not meet at $z$ in
$\g$, we have that the entry $r$ must have mask-value 1.  Hence, the only way
for the string which is in the third position at the top of $\g$ to end up on
the bottom of $\g$ in a position lying to the left of the strings labeled 2 and
3 is if $\g$ has a plain-zero $y$ as shown in Figure~\ref{f:mu_k3cd}.  The
position of the plain-zero $y$ in $\g$ is completely determined by the position
of the string labeled 1 in $\s$ since $w^{\s} = w^{\g}$.  By
Lemma~\ref{l:no_isolation} there exists a zero-defect $e$ with $y$ as its left
critical zero.  If $e$ lies further from $y$ than shown, we obtain an I-shape in
$\g$ after taking account of the necessary string crossings for $d$ and $e$.
Hence, the right string of $d$ becomes the left string of $e$, and the position
of $t$ in $\g$ is completely determined by the position of the string labeled 2
in $\s$.  Also, the right critical zero $v$ of $e$ occurs in a column weakly
left of the column of $t$, for otherwise the strings of $e$ will not have an
opportunity to cross again.

Observe that $y$, $t$ and $v$ must be plain-zeros since if any of these entries
were zero-defects, the string crossing would introduce a 4-stack into the heap,
contradicting that $w$ is Deodhar.  By Lemma~\ref{l:no_isolation}, there exists
a zero-defect $f$ with $t$ as the left critical zero.  Using Lemma~\ref{l:d.pi},
we can project $\g$ to type $A$ and obtain a $\mu$-mask on a Deodhar
permutation.  Moreover, the projection does not alter the mask-values of the
entries to the right of column 1.  Consequently, if the zero-defect $e$ does not
share a critical zero with another zero-defect to the right, then $e$ and $f$
form a separating pair in the projection of $\g$, which is a contradiction.  But
then we obtain an I-shape in the heap of $w$, contradicting that $w$ is Deodhar.

Finally, it remains to show that if $\s$ is in case (c) or (d) then $\g$ cannot
fall into case (a).  If $\g$ is in case (a) then $\g$ has a decorated heap 
of the form:
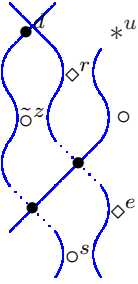
\[
\heap {
& & \StringLRX{\hf^d} & & {\hv^u} & & {\hs} & \\
& \StringR{\hs} & & \StringLR{\hd^r} & & {\hs} & & {\hs} \\
{\hs} & & \StringLR{\tilde{\hz}^z} & & \StringL{\hz} & & {\hs} & \\
& \StringR{\hs} & & \StringLXDR{\hf} & & {\hs} & & {\hs} \\
& & \StringLXDR{\hf} & & \StringL{\hd^e} & & {\hs} & \\
& {\hs} & & \StringLR{\hz^s} & & {\hs} & & {\hs} \\
}
\]

Since $w^{\s} = w^{\g}$ we have that there are no negatively labeled strings in
$\g$ so $z$ must be a plain-zero.  Also, the labeled strings at $d$ must form an
inversion in the 1-line notation for $x = w^{\s}$.  Thus, the path of the left
string of $d$ must encounter a plain-zero $s$ along the minimal southwest
diagonal in order for the strings of $d$ to be labeled in $\g$ as they are in
$\s$.  Since $s$ is a plain-zero, we have by Lemma~\ref{l:no_isolation} that
there exists a zero-defect $e$ with $s$ as its left critical zero.  Hence, $u$
cannot be in the heap as it would form an I-shape with the string crossing for
$e$ by lateral convexity.  Since $z$ is a plain-zero, we have by
Lemma~\ref{l:no_isolation} that $r$ must be a zero-defect with $z$ as the left
critical zero.  

Because $d$ is a zero-defect in $\s$, we have that the labeled strings lying in
the first and second position at the top of the decorated heap of $\s$ form an
inversion in the 1-line notation for $x = w^{\s}$.  Since $r$ is zero-defect in
$\g$, the labeled strings lying in the first and third position at the top of
the decorated heap of $\g$ form an inversion in the 1-line notation for $x =
w^{\g}$.  Hence, the labeled strings in the first three positions at the top of
the decorated heaps of both $\s$ and $\g$ either form a $[321]$ permutation
pattern, or a $[312]$ permutation pattern.

If the labeled strings in the first three positions at the top of the decorated
heaps form a $[321]$ pattern then the decorated heap of $\g$ has the form shown
in Figure~\ref{f:mu_k3gcd}(1).  In particular, the string labeled 1 ends up left
of $s$ on the bottom of $\g$.  Also, neither of the strings labeled 2 or 3 can
encounter a mask-value 0 entry below $z$ until after crossing the string labeled
1 in $\g$.  Hence, the right critical zero of $r$ must be a plain-zero.  Using
Lemma~\ref{l:d.pi}, we can project $\g$ to type $A$ and obtain a $\mu$-mask on a
Deodhar permutation.  Moreover, the projection does not alter the mask-values of
the entries to the right of column 1.  Consequently, if the zero-defect $r$ does
not share its right critical zero with another zero-defect, then $e$ and $r$
form a separating pair in the projection of $\g$, which is a contradiction.  But
then we obtain an I-shape in the heap of $w$, contradicting that $w$ is Deodhar.

\begin{figure}[h]
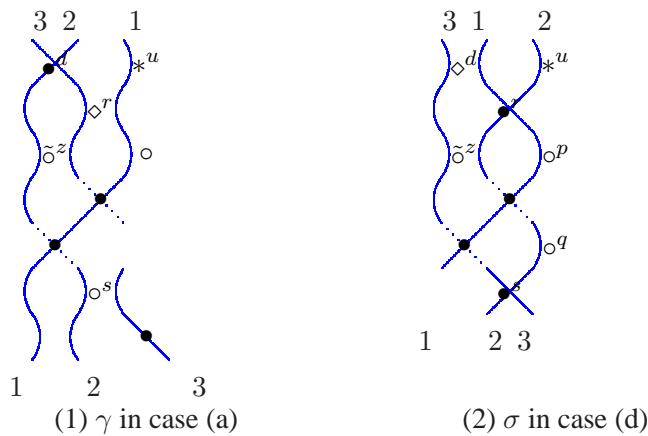

\begin{tabular}{cc}
\heap {
& {\hs} & 3 \ \ 2 & & 1 \ \ \  & {\hs} & & {\hs} \\
& & \StringLRX{\hf^d} & & \StringL{\hv^u} & & {\hs} & \\
& \StringR{\hs} & & \StringLR{\hd^r} & & {\hs} & & {\hs} \\
{\hs} & & \StringLR{\tilde{\hz}^z} & & \StringL{\hz} & & {\hs} & \\
& \StringR{\hs} & & \StringLXDR{\hf} & & {\hs} & & {\hs} \\
& & \StringLXDR{\hf} & & {\hs} & & {\hs} & \\
& \StringR{\hs} & & \StringLR{\hz^s} & & {\hs} & & {\hs} \\
& & \StringLR{\hs} & & \StringLX{\hf} & & {\hs} & \\
& \ \ 1  & & 2 \ \ & &  \ \ 3 & {\hs} & & {\hs} \\
} &
\heap {
& {\hs} & 3 \ \ 1 & & 2 \ \ \  & {\hs} & & {\hs} \\
& & \StringLR{\hd^d} & & \StringL{\hv^u} & & {\hs} & \\
& \StringR{\hs} & & \StringLRX{\hf^r} & & {\hs} & & {\hs} \\
{\hs} & & \StringLR{\tilde{\hz}^z} & & \StringL{\hz^p} & & {\hs} & \\
& \StringR{\hs} & & \StringLXDR{\hf} & & {\hs} & & {\hs} \\
& & \StringLXDR{\hf} & & \StringL{\hz^q} & & {\hs} & \\
& {\hs} & & \StringLRX{\hf^s} & & {\hs} & & {\hs} \\
& \ \ 1  & & 2 \ \ 3 & {\hs} & & {\hs} & & {\hs} \\
} \\
(1) $\g$ in case (a) & (2) $\s$ in case (d) \\
\end{tabular}
\caption{$\g$ in case (a) distinguished from $\s$ in case (c) or (d)}
\label{f:mu_k3gcd}
\end{figure}

Suppose the labeled strings in the first three positions at the top of the
decorated heaps form a $[312]$ pattern.  By the remarks justifying
Figure~\ref{f:mu_k3cd}, we find that $\s$ cannot be in case (c) so $\s$ must be
in case (d).  Hence, the decorated heap of $\s$ has the form shown in
Figure~\ref{f:mu_k3gcd}(2).  The string labeled 3 cannot encounter a mask-value
0 entry below $z$ in $\s$ until it crosses the strings labeled 1 and 2.  Hence,
the position of $q$ in $\s$ is completely determined by the position of the
string labeled 2 on the bottom of $\g$.  Thus, the position of the string
labeled 1 in $\g$ determines the location of the right critical zero $p$ of $d$
in $\s$ and $p$ occurs in a column weakly left of the column of $q$.  Since $q$
must be a plain-zero, we have that there exists a zero-defect $f$ northeast of
$q$ by Lemma~\ref{l:no_isolation}.  Using Lemma~\ref{l:d.pi}, we can project
$\s$ to type $A$ and obtain a $\mu$-mask on a Deodhar permutation.  Moreover,
the projection does not alter the mask-values of the entries to the right of
column 1.  Consequently, if the zero-defect $d$ does not share its plain right
critical zero $p$ with another zero-defect, then $d$ and $f$ form a separating
pair in the projection of $\s$, which is a contradiction.  But then we obtain an
I-shape in the heap of $w$, contradicting that $w$ is Deodhar.

Thus, we have shown that $\s$ and $\g$ always fall into the same case (a) - (d).
If $\s$ and $\g$ fall into case (a) then we can remove the maximal entry $d$
which is a plain-one in both masks, and then apply the argument from the case
$k(w)=2$ to obtain a contradiction.  

Next, assume that $\s$ and $\g$ fall into case (b), (c) or (d).  We apply
Lemma~\ref{l:d.pi} to add three new entries to the decorated heaps of $\s$ and
$\g$, obtaining a pair of decorated heaps $\check{\s}$ and $\check{\g}$ on the
Deodhar type $A$ element $\pi(w)$.  The mask-values of the three new entries
depend on which case (b), (c) or (d) the masks fall into, but in each case the
projection satisfies $D(\check{\s}) = D(\check{\g})$.  To verify that
$\pi(w)^{\check{\s}} = \pi(w)^{\check{\g}}$ we must either show that the same
strings encounter the three new entries from the top of the decorated heap of
$\s$ as in $\g$, or we must show that the same strings encounter the three new
entries from the bottom of the decorated heap of $\s$ as in $\g$.

Consider the heap of $w$ with entries labeled as follows.
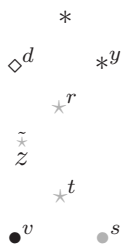
\[
\xymatrix @=-4pt @! {
& {\hs} & & {\hv} & & {\hs} & & {\hs} \\
& & {\hd^d} & & {\hv^y} & & {\hs} & \\
& {\hs} & & {\ha^r} & & {\hs} & & {\hs} \\
{\hs} & & {\stackrel{\tilde{\ha}}{z}} & & {\hs} & & {\hs} & \\
& {\hs} & & {\ha^t} & & {\hs} & & {\hs} \\
& & {\hf^v} & & {\hb^s} & & {\hs} & \\
& {\hs} & & {\hs} & & {\hs} & & {\hs} \\
} 
\]
If both $s$ and $y$ lie in the heap of $w$ then we introduce an I-shape,
contradicting that $w$ is Deodhar.  Hence, at least one of $\{ s, y \}$ is not
in the heap of $w$.  If $s$ is not in the heap of $w$ then the mask-value of $t$
in $\s$ and $\g$ is completely determined by the case (b) - (d) that $\s$ and
$\g$ fall into and the requirement that the strings of $d$ cross at $v$.  In
particular, the strings encountering $v$, $t$, and $z$ from the bottom of the
decorated heap of $\s$ must be the same as the strings encountering $v$, $t$,
and $z$ from the bottom of the decorated heap of $\g$.  Hence,
$\pi(w)^{\check{\s}} = \pi(w)^{\check{\g}}$ in this case.

On the other hand, if $y$ is not in the heap of $w$ then we may observe that the
mask values of $d$, $z$ and $r$ in $\s$ and $\g$ are completely determined by
the case (b) - (d) that $\s$ and $\g$ fall into.  Thus, we have that the strings
encountering $d$, $r$ and $z$ from the top of the decorated heap of $\s$ must be
the same as the strings encountering $d$, $r$ and $z$ from the top of the
decorated heap of $\g$.  Hence, $\pi(w)^{\check{\s}} = \pi(w)^{\check{\g}}$ in
this case.

Since $\s$ and $\g$ are distinct, we obtain a contradiction to
Theorem~\ref{t:mu.main}.

\bigskip

\noindent \textbf{Case $\mathbf{k(w)=4}$}.  Suppose that $w$ has entries on 4
distinct levels in the first column, and consider the decorated heap on $w$ that 
is associated to the $\mu$-mask $\s$.  Following the argument in the proof of
\cite[Proposition 8.7]{b-j1}, we reduce to the case where $\s$ contains a
decorated heap fragment of the form:

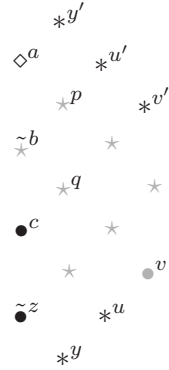
\begin{equation}\label{e:pi.map.4}
\xymatrix @=-4pt @! {
& {\hv}^{y'} & & {\hs} & & {\hs} & & {\hs} \\
{\hd}^a & & {\hv}^{u'} & & {\hs} & & {\hs} \\
& {\ha}^p & & {\hv}^{v'} & & {\hs} & & {\hs} \\
{\tilde{\ha}}^b & & {\ha} & & {\hs} & & {\hs} \\
& {\ha}^q & & {\ha} & & {\hs} & & {\hs} \\
{\hf}^c & & {\ha} & & {\hs} & & {\hs} \\
& {\ha} & & {\hb}^{v} & & {\hs} & & {\hs} \\
{\tilde{\hf}}^z & & {\hv}^{u} & & {\hs} & & {\hs} \\
& {\hv}^{y} & & {\hs} & & {\hs} & & {\hs} \\
} 
\end{equation}

If there is a point at $y$ or $y'$ then there is a 4-stack in the second
column.  If there is a point at $u$ or $u'$, we obtain an I-shape with the other
entries that exist by convexity.  In either case, we contradict that $w$ is
Deodhar.  Similarly, if both $v$ and $v'$ exist, then we obtain an I-shape, so
at least one of them is not in the heap.  By considering $w^{-1}$ if necessary,
we can assume that $v'$ is not in the heap.

By Lemma~\ref{l:no_isolation}, we can assume that $z$ has mask-value 1 and that
it is a string crossing for some zero-defect, since if $z$ has mask-value 0 then
it can play no critical role in any defect.  Also by Lemma~\ref{l:no_isolation},
we have that $a$ is never a plain-zero in any $\mu$-mask.  Thus, if $a$ is not a
zero-defect in at least one of $\{ \s, \g \}$ then it has mask-value 1 in both
masks so we can remove it from both masks and obtain a contradiction by
considering the $k(w)=3$ case.  Therefore, we assume without loss of generality
that $a$ is a zero-defect in $\s$.  This forces $c$ to have mask-value 1 since
$a$ cannot be a defect otherwise.  Thus, we have a heap fragment of the form in
\eqref{e:pi.map.4}.

Also, if $b$ is a plain-zero then $z$ cannot be the string crossing of any
zero-defect, regardless of the mask-value for $a$.  Hence, we may assume that
$b$ is either a zero-defect or has mask-value 1 in all $\mu$-masks on $w$.

\begin{figure}[h]
\begin{tabular}{cccc}
\heap {
 & \bar{1} \ \ \bar{2} & \\
& \StringLR{\hd^a} & & \\
\StringR{\hs} & & \StringL{\hz^p} \\
& \StringLR{\tilde{\hd}^b} \\
\StringR{\hs} & & \StringL{\hz} & \\
& \StringLRX{\stackrel{c}{\hf}} & \\
\StringR{\hs} & & \StringL{\hz} & \\
& \StringLRX{\stackrel{z}{\tilde{\hf}}} & \\
 & 1 \ \ 2 & {\hs} & \\
} &
\heap {
& 2 \ \ \bar{1} &  & \\
& \StringLR{{\hd}^a} & & {\hs} & \\
\StringR{\hs} & & \StringLX{\hf^p} & & {\hs} & \\
& \StringLX{\tilde{\hf}^b} & & \StringLX{\hf} & \\
{\hs} & & \StringL{\hz^q} & & \StringL{\hz} & \\
& \StringRX{\hf^c} & & \StringRX{\hf} & \\
\StringR{\hs} & & \StringRX{\hf} & & {\hs} & \\
& \StringLRX{\tilde{\hf}^z} & & {\hs} & \\
& 1 \ \ 2 &  & & {\hs} & \\
} &
\heap {
& 1 \ \ \bar{2} & & \\
& \StringLR{{\hd}^a} & \\
\StringR{\hs} & & \StringLX{\hf^p} & & {\hs} & \\
& \StringLX{\tilde{\hf}^b} & & \StringL{\hz} & \\
{\hs} & & \StringLRX{\hf^q} & & {\hb} & \\
& \StringRX{\hf^c} & & \StringL{\hz} & \\
\StringR{\hs} & & \StringRX{\hf} & & {\hs} & \\
& \StringLRX{\tilde{\hf}^z} & \\
& 1 \ \ 2 & & & {\hs} & \\
} &
\heap {
& 2 \ \ 1 & & & {\hs} & \\
& \StringLR{{\hd}^a} & \\
\StringR{\hs} & & \StringL{\hz^p} & \\
& \StringLRX{\tilde{\hf}^b} & & {\hb} & \\
\StringR{\hs} & & \StringL{\hz^q} & & {\hb} \\
& \StringLRX{\hf^c} & & {\hb} \\
\StringR{\hs} & & \StringL{\hz} & & {\hs} \\
& \StringLRX{\tilde{\hf}^z} \\
& 1 \ \ 2 & & & {\hs} \\
& {\hs} & & {\hs} \\
} \\
(a) & (b) & (c) & (d) \\
\end{tabular}
\caption{$\s$ is a $\mu$-mask with four levels in column 1}
\label{f:mu_k4}
\end{figure}

Working through all of the cases given in the proof of \cite[Proposition
8.7]{b-j1} and applying Lemma~\ref{l:no_isolation}, we find that $\s$ must have
one of the decorated heap fragments shown in Figure~\ref{f:mu_k4}.  Note that
since the labeled strings are distinct among these four cases and $w^{\s} = x =
w^{\g}$, we have that $\s$ falls into one of these cases if and only if $\g$
falls into the same case.

\begin{figure}[h]
\begin{tabular}{cccc}
\heap {
 & \bar{1} \ \ \bar{2} & & \\
& \StringLRX{\hf^a} & \\
\StringR{\hs} & & \StringL{\hd^p} & \\
& \StringLR{\tilde{\hd}^b} & & \StringL{\hz} \\
\StringR{\hs} & & \StringLRX{\hf^q} & \\
& \StringL{\stackrel{c}{\hz}} & \\
\StringR{\hs} & & \StringRX{\hf} & \\
& \StringLRX{\stackrel{z}{\tilde{\hf}}} & \\
 & 1 \ \ 2 & {\hs} & \\
} &
\heap {
& 2 \ \ \bar{1} &  & \\
& \StringLRX{{\hf}^a} & & {\hs} & \\
\StringR{\hs} & & {\hs} & & {\hs} & \\
& \StringL{\tilde{\hz}^b} & & {\hs} & \\
\StringR{\hs} & & {\hs} & & {\hs} & \\
& \StringLX{\hf^c} & & {\hs} & \\
{\hs} & & \StringL{\hz} & & {\hs} & \\
& \StringLRX{\tilde{\hf}^z} & & {\hs} & \\
& 1 \ \ 2 &  & & {\hs} \\
} &
\heap {
& 1 \ \ \bar{2} & & \\
& \StringLRX{{\hf}^a} & \\
\StringR{\hs} & & \StringL{\hd^p} & \\
& \StringLR{\tilde{\hz}^b} & \\
\StringR{\hs} & & {\hs} & \\
& \StringL{\hz^c} & \\
\StringR{\hs} & & {\hs} \\
& \StringLRX{\tilde{\hf}^z} & \\
& 1 \ \ 2 & & \\
} &
\heap {
& 2 \ \ 1 & \\
& \StringLRX{{\hf}^a} & \\
\StringR{\hs} & & \StringL{\hz^p} \\
& \StringLRX{\tilde{\hf}^b} \\
\StringR{\hs} & & {\hs} \\
& \StringL{\hz^c} \\
\StringR{\hs} \\
& \StringLRX{\tilde{\hf}^z} \\
& 1 \ \ 2 & & \\
& {\hs} & \\
} \\
(a) & (b) & (c) & (d) \\
\end{tabular}
\caption{$\g$ is a $\mu$-mask with four levels in column 1}
\label{f:mu_k4g}
\end{figure}

Next, suppose that $\g$ has $a$ with mask-value 1 and $\s$ falls into one of the
cases $\{$ (a), (b), (c), (d) $\}$.  Then we claim $\g$ has one of the forms
shown in Figure~\ref{f:mu_k4g}, corresponding to the case that $\s$ falls into.
In any $\mu$-mask we have that $z$ is a plain-one.  Since $w^{\s} = x = w^{\g}$,
the string paths in $\g$ determine the mask-values of the entries in column 1.
In case (a), we have $b$ with mask-value 0 in order that the strings of $a$ be
labeled negatively.  Hence, $b$ must be a zero-defect whose strings cross at
$z$.  This implies that the strings of $a$ become the strings of $b$, so by
Lemma~\ref{l:no_isolation} there exists a zero-defect $p$ as shown.  But this is
a contradiction because $p$ is a plain-zero in $\s$ and $s_2$ is an ascent for
the element $x s_1$.

In case (b), the string labeled 1 forces $b$ to have mask-value 0, but then the
other string of $a$ must be labeled negatively, contradicting that $w^{\s} = x
= w^{\g}$.  Similarly, in case (c), the string labeled 2 forces $b$ to have
mask-value 0, but then the other string of $a$ must be labeled negatively,
contradicting that $w^{\s} = x = w^{\g}$.  In case (d), note that $p$ cannot be
a defect because the first and last opportunity for the strings of $p$ to cross
occurs at $z$ which contradicts that the strings through $a$ are assumed to
cross at $z$.  But then $p$ is a plain-zero which contradicts
Lemma~\ref{l:no_isolation}.

Thus, we have shown that $\s$ and $\g$ always fall into the same case.  As in
the proof of \cite[Proposition 8.7]{b-j1}, we may remove the zero-defect $a$ and
change the mask-value of some plain-zero entry to 1.  Specifically, in cases
(a), (c) and (d), we change the mask-value of the right critical zero for $a$.
Since the entry directly northeast of the right critical zero of $a$ is not in
the heap for any of these cases, the same strings encounter the entry from the
top of the decorated heap of $\s$ as in $\g$.  In case (b), we change the
mask-value of $q$ to have mask-value 1.  Since there is no entry two units
northeast from $q$ and the mask-value of the entry directly northeast of $q$ is
determined, the same strings encounter $q$ from the top of the decorated heap of
$\s$ as in $\g$.

We obtain in this way a pair of decorated heaps $\check{\s}$ and $\check{\g}$ on
a Deodhar element $\check{w}$ with $k(w) = 3$.  These reductions remove a single
zero-defect and a single plain-zero in each case, so $D(\check{\s}) =
D(\check{\g})$.  Because $w^{\s} = w^{\g}$ and the same strings encounter the
entries that we changed in both masks we have $\check{w}^{\check{\s}} =
\check{w}^{\check{\g}}$.  Thus, we may apply the argument from the case $k(w)=3$
to obtain a contradiction to Theorem~\ref{t:mu.main}.

\bigskip

\noindent \textbf{Case $\mathbf{k(w)=5}$}.  If $w$ has entries
on 5 or more distinct levels in the first column then $w$ is not Deodhar by
Lemma~\ref{l:shape}.

\bigskip

This exhausts the cases, completing the proof of the proposition.
\end{proof}

\begin{proposition}\label{p:mu.d.not.convex}
Suppose that $w \in D_n$ is a Deodhar element which is not convex.  Then,
$\mu(x, w) \in \{0, 1\}$ for all $x < w$.
\end{proposition}
\begin{proof}
Since $w$ is not convex, there exists a minimal pair of entries in column $i
\geq 2$ with only a left resolution.  By \cite[Lemma 6.8]{b-w}, the part of the
heap to the left of column $i$ has the particular form shown in
Figure~\ref{f:mu_not_convex} in which each column $1, \ldots, i$ has exactly two
entries.  

\begin{figure}[h]
\begin{tabular}{c}
\xymatrix @=-8pt @! {
& {\hs} & & {\hf^x} & & {\hs} & & {\hs} \\
{\hs} & & {\hf} & & {\hs} & & {\hs} \\
& {\hf} & & {\hs} & & {\hs} & & {\hs} \\
{\hf \tilde{\hf}} & & {\hs} & & {\hs} & & {\hs} \\
& {\hf} & & {\hs} & & {\hs} & & {\hs} \\
{\hs} & & {\hf} & & {\hs} & & {\hs} \\
& {\hs} & & {\hf^y} & & {\hs} & & {\hs} \\
}
\end{tabular}
\caption{Fragment from a non-convex element}
\label{f:mu_not_convex}
\end{figure}
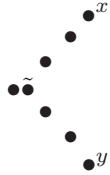

Suppose for the sake of contradiction that there exist distinct $\mu$-masks $\s$
and $\g$ such that $w^{\s} = x = w^{\g}$.  If exactly one of the entries in
column 1 has mask-value 0 in $\s$ or $\g$, then the mask-value 0 entry cannot be
the left critical zero for any zero-defect because the left string of such a
zero-defect will be negative while the right string will be positive, and this
would contradict Lemma~\ref{l:no_isolation}.  Hence, both entries in column 1
must have the same mask-value in any $\mu$-mask.  Then, both entries of column 1
have mask-value 0 in $\s$ if and only if all of the entries in the 1-line
notation for $x = w^{\s}$ are positive, and this occurs if and only if both
entries of column 1 have mask-value 0 in $\g$ because $w^{\s} = x = w^{\g}$.
Thus, the mask-values of $\s$ and $\g$ must agree on the entries in column 1.  

If the lower entry of column $j<i$ in the heap fragment has mask-value 0, then
since it cannot be a defect, it must enable the string crossing for some
zero-defect lying to the right of column $j$ by Lemma~\ref{l:no_isolation}.  In
particular, by the argument given in the proof of \cite[Lemma 6.10]{b-j1}, the
top entry of column $j$ must also be a plain-zero.  Hence, the mask-values of
$\s$ and $\g$ agree on all columns in which the lower entry has mask-value 0 for
one of the masks.

Next, observe that each column in the heap fragment in which both entries have
mask-value 1 corresponds to a fixed point in the 1-line notation for $x$,
viewing $x$ as a signed permutation.  Therefore, column $j$ has both entries
with mask-value 1 in $\s$ if and only if column $j$ has both entries with
mask-value 1 in $\g$, for each $1 \leq j \leq i$.

Thus, the mask-values of $\s$ and $\g$ must agree in every column $1, \ldots, i$.
The proof of \cite[Lemma 6.10]{b-w} describes how to project any mask on $w$ to
a mask $\check{\s}$ on a convex element $\check{w}$ obtained from $w$ by
removing columns $1, \ldots, i$.  This projection has the properties that:
\begin{enumerate}
\item[(1)]  $D(\s) = D(\check{\s})$.
\item[(2)]  The mask-values of the entries to the right of column $i$ in the
	decorated heap of $\s$ are unchanged in the decorated heap of $\check{\s}$.
\end{enumerate}
Since $w^{\s} = w^{\g}$ we also have that $\check{w}^{\check{\s}} =
\check{w}^{\check{\g}}$.  Because $\s$ and $\g$ are distinct, we obtain a
contradiction to Proposition~\ref{p:mu.d.convex}.
\end{proof}

Combining Propositions~\ref{p:mu.d.convex} and \ref{p:mu.d.not.convex}, we have
the following result.

\begin{theorem}\label{t:main.d}
Let $w$ be a Deodhar element of type $D$ and suppose $x$ is another type $D$ element.  Then,
$\mu(x,w) \in \{0, 1\}$.
\end{theorem}


\bigskip
\section{The 0-1 result for Deodhar elements of finite exceptional types}\label{s:mu.others}

\begin{theorem}\label{t:main.e}
If $w$ is a Deodhar in type $E_6$, $E_7$, $E_8$, $F_4$ or $G_2$ then $\mu(x,w)
\in \{0,1\}$ for all $x$.
\end{theorem}
\begin{proof}
These groups are finite, so the computation of $\mu$ values is verifiable by
computer.   Some C++ code is available at
\url{http://www.math.ucdavis.edu/~brant/code/} for this purpose.
\end{proof}

This completes the proof of Theorem~\ref{t:mu.remain}.


\bigskip
\section*{Acknowledgments}

We thank Sara Billey, Richard Green, Monty McGovern, Yuval Roichman, Monica
Vazirani and Greg Warrington for many useful conversations and suggestions.


\appendix
\section{}\label{s:appendix}

Here is an explicit description of the map defined in Lemma~\ref{l:d.pi} that
takes a decorated heap of type $D$ to a decorated heap of type $A$.  The map
adds three new entries to the left of the first column of the type $D$ heap so
that the result is a fully-commutative heap in type $A$.  The table below
indicates the mask-values that should be assigned to these new entries and the
entries from the first column of the type $D$ heap so that the Deodhar statistic
is preserved.  This assignment of mask-values depends on the cases given in the
first column.  Further details are given in \cite[Proposition 8.7]{b-j1}.

\bigskip
\begin{tabular}{|p{1.8in}|c|c|}
	\hline
Case & $(w, \s)$ & $(\pi(w), \pi(\s))$ \\
	\hline
$d$ is a plain-one. & 
\heap{
& {\hs} & & \StringL{\hv} & & {\hs} & & {\hs} \\
& & \StringRX{\stackrel{\hf}{d}} & & {\hs} & & {\hs} & \\
& {\hs} & & {\ha} & & {\hs} & & {\hs} \\
{\hs} & & \tilde{\ha} &  & {\ha} & & {\hs} & \\
& {\hs} & & {\ha} & & {\hs} & & {\hs} \\
& & {\ha} & & {\hs} & & {\hs} & \\
} &
\heap {
& {\hs} & & \StringL{\hv} & & {\hs} & & {\hs} \\
& & \StringRX{\stackrel{\hf}{d}} & & {\hs} & & {\hs} & \\
& \StringRX{\hf} & & {\ha} & & {\hs} & & {\hs} \\
\StringRX{\hf} & & \tilde{\ha} &  & {\ha} & & {\hs} & \\
& {\hf} & & {\ha} & & {\hs} & & {\hs} \\
& & {\ha} & & {\hs} & & {\hs} & \\
} \\ 
	\hline
$d$ is a zero-defect and $z$ has mask-value 1. &
\heap {
& {\hs} & & {\hv} & & {\hs} & & {\hs} \\
& & \StringLR{\hd^d} & & {\hs} & & {\hs} & \\
& \StringR{\hs} & & \StringL{\hz} & & {\hs} & & {\hs} \\
{\hs} & & \StringLRX{\stackrel{\tilde{\hf}}{z}} & & {\hs} & & {\hs} & \\
& \StringR{\hs} & & \StringLXD{\hf} & & {\hs} & & {\hs} \\
& & \StringLX{\hf} & & \StringL{\hz} & & {\hs} & \\
& {\hs} & & \StringLRX{\hf} & & {\hs} & & {\hs} \\
} &
\heap {
& {\hs} & & {\hv} & & {\hs} & & {\hs} \\
& & \StringLR{\hd^d} & & {\hs} & & {\hs} & \\
& \StringLR{\hd} & & \StringL{\hz} & & {\hs} & & {\hs} \\
\StringR{\hz} & & \StringLRX{\stackrel{\hf}{z}} & & {\hs} & & {\hs} & \\
& \StringLRX{\hf} & & \StringLXD{\hf} & & {\hs} & & {\hs} \\
& & \StringLX{\hf} & & \StringL{\hz} & & {\hs} & \\
& {\hs} & & \StringLRX{\hf} & & {\hs} & & {\hs} \\
} \\
	\hline
$d$ is a zero-defect and $z$ has mask-value 0.  The strings for $d$
meet at $z$. &
\heap {
& {\hs} & & {\hv} & & {\hs} & & {\hs} \\
& & \StringLR{\hd^d} & & {\hs} & & {\hs} & \\
& \StringR{\hs} & & \StringL{\hz} & & {\hs} & & {\hs} \\
{\hs} & & \StringLR{\stackrel{\tilde{\hz}}{z}} & & {\hs} & & {\hs} & \\
& \StringR{\hs} & & {\hb} & & {\hs} & & {\hs} \\
& & \StringLX{\hf} & & {\hs} & & {\hs} & \\
} &
\heap {
& {\hs} & & {\hv} & & {\hs} & & {\hs} \\
& & \StringLR{\hd^d} & & {\hs} & & {\hs} & \\
& \StringR{\hz} & & \StringL{\hz} & & {\hs} & & {\hs} \\
{\hf} & & \StringLR{\stackrel{\hd}{z}} & & {\hs} & & {\hs} & \\
& \StringR{\hz} & & {\hb} & & {\hs} & & {\hs} \\
& & \StringLX{\hf} & & {\hs} & & {\hs} & \\
} \\
	\hline
\end{tabular}

\begin{tabular}{|p{1.8in}|c|c|}
	\hline
Case & $(w, \s)$ & $(\pi(w), \pi(\s))$ \\
	\hline
$d$ is a zero-defect and $z$ has mask-value 0.  The strings for $d$
do not meet at $z$. &
\heap {
& {\hs} & & {\hv} & & {\hs} & & {\hs} \\
& & \StringLR{\hd^d} & & {\hs} & & {\hs} & \\
& \StringR{\hs} & & \StringLX{\hf} & & {\hs} & & {\hs} \\
{\hs} & & \StringL{\stackrel{\tilde{\hz}}{z}} & & {\hs} & & {\hs} & \\
& \StringR{\hs} & & {\hb} & & {\hs} & & {\hs} \\
& & \StringLX{\hf} & & {\hs} & & {\hs} & \\
} &
\heap {
& {\hs} & & {\hv} & & {\hs} & & {\hs} \\
& & \StringLR{\hd^d} & & {\hs} & & {\hs} & \\
& \StringRX{\hf} & & \StringLX{\hf} & & {\hs} & & {\hs} \\
\StringR{\hz} & & {\stackrel{\hf}{z}} & & {\hs} & & {\hs} & \\
& \StringLX{\hf} & & {\hb} & & {\hs} & & {\hs} \\
& & \StringLX{\hf} & & {\hs} & & {\hs} & \\
} \\
	\hline
\end{tabular}


\bibliographystyle{alpha}
\bibliography{our}

\end{document}